\documentclass[preprint,review,3p,12pt]{elsarticle}

\usepackage{amssymb,dsfont,amsmath}
\usepackage{lineno}
\usepackage{graphicx}
\usepackage{xcolor}
\usepackage{epstopdf}
\usepackage{bm}

\newtheorem{definition}{Definition}
\newtheorem{principle}{Principle}
\biboptions{comma,sort}

\everymath{\displaystyle}

\journal{Journal of Computational Physics}

\def\bigdot{\bm{\cdot}}
\DeclareMathOperator*{\arginf}{arg\,inf}

\begin{document}

\begin{frontmatter}
  \title{A least-squares implicit RBF-FD closest point method and applications to PDEs on moving surfaces}
  \author[rvt,rvt4]{A.~Petras\corref{cor1}}
  \ead{apetras@bcamath.org}
  \author[rvt3]{L.~Ling}
  \ead{lling@hkbu.edu.hk}
  \author[rvt2]{C.~Piret}
  \ead{cmpiret@mtu.edu}
  \author[rvt]{S.J.~Ruuth}
  \ead{sruuth@sfu.ca}
  \cortext[cor1]{Corresponding author}
  \address[rvt]{Department of Mathematics, Simon Fraser University, Burnaby, British Columbia, Canada V5A1S6}
  \address[rvt2]{Department of Mathematical Sciences, Michigan Technological University, Michigan, USA}
  \address[rvt3]{Department of Mathematics, Hong Kong Baptist University, Kowloon Tong, Hong Kong}
  \address[rvt4]{BCAM-Basque Center for Applied Mathematics, Bilbao, Basque Country, Spain 48009}

  \begin{abstract}
The closest point method (Ruuth and Merriman, J. Comput. Phys. 227(3):1943–-1961, [2008]) is an embedding method
developed to solve a variety of partial differential equations (PDEs) on smooth surfaces, using a closest point representation of the surface and
standard Cartesian grid methods in the embedding space.
Recently, a closest point method with explicit time-stepping was proposed that uses finite differences derived from radial basis functions (RBF-FD).
Here, we propose a least-squares implicit formulation of the closest point method to impose the constant-along-normal extension
of the solution on the surface into the embedding space.
Our proposed method is particularly flexible with respect to the choice of the
computational grid in the embedding space.
In particular, we may compute over a computational tube that contains problematic nodes.
This fact enables us to combine the proposed method with the grid based particle method
(Leung and Zhao, J. Comput. Phys. 228(8):2993--3024, [2009]) to obtain a numerical method for approximating PDEs on moving surfaces.
We present a number of examples to illustrate the numerical convergence properties of our proposed method.
Experiments for advection-diffusion equations and Cahn-Hilliard equations that are strongly coupled to the velocity of the surface are also presented.
\end{abstract}

\end{frontmatter}

\section{Introduction}
Many applications in the natural and applied sciences involve the solution of partial differential equations (PDEs) on surfaces.
Application areas for PDEs on static surfaces include image processing \cite{bertalmio2001navier,tian2009segmentation,biddle2013volume},
biology \cite{olsen1998spatially,murray2001mathematical} and computer graphics \cite{auer2012real}.
Applications for
PDEs on moving surfaces also occur frequently.
Notable examples arise in biology \cite{elliott2010modeling,elliott2012modelling,barreira2011surface,venkataraman2011modeling},
material science \cite{eilks2008numerical}, fluid dynamics \cite{adalsteinsson2003transport,james2004surfactant} and computer
graphics \cite{auer2013semi}.

Methods for solving PDEs on surfaces can be categorized according to the representation of the surface. On static surfaces, there are methods that
solve PDEs on parametrized surfaces \cite{lui2005solving,floater2005surface}, and on triangulated surfaces using finite difference
\cite{turk1991generating} or finite element \cite{dziuk2007surface} methods.  Also popular are the embedding methods,
which solve PDEs on surfaces embedded in a higher dimensional space using a projection operator \cite{greer2006improvement,flyer2009radial,piret2012orthogonal,fuselier2013high}.

Similarly, methods for solving PDEs on moving surfaces can be categorized according to the representation of the surface.
On moving triangular meshes, some commonly used methods for solving PDEs are the finite element method \cite{dziuk2007finite,dziuk2013finite}
and the finite volume method \cite{nemadjieu2012finite}.  On parametrized surfaces and surfaces represented by particles, a direct
discretization of the parametrized differential operators can be applied \cite{leung2011grid,leung2010gaussian}. Finally, on surfaces embedded
in higher dimensional spaces, the zero level set of a function defined in the embedding space is commonly used to represent surfaces
\cite{sethian1999level,osher2006level}. Methods for solving PDEs on such surfaces include finite element methods \cite{dziuk2010eulerian} and
finite differences \cite{xu2003eulerian}.

For surfaces embedded in a higher dimensional space,
an alternative way to represent the surface is to compute and store the closest points to the surface over a neighborhood of the surface.
The closest point method (CPM) \cite{ruuth2008simple} is an embedding method for solving PDEs on surfaces that uses such a representation. The surface
differential operators are replaced with Cartesian ones by extending the solution to the embedding space using interpolation and the closest point
representation. Then, standard finite difference schemes are used to solve the PDE on the surface. A closest
point method with explicit temporal discretization is also available \cite{macdonald2009implicit}, allowing the use of large time step-sizes in implicit time-stepping discretizations.
Note also that  the closest point method may be combined with a modified grid based particle method  to
solve PDEs on {\it moving} surfaces; see \cite{petras2016pdes} for details.

Recently, a closest point method using finite differences derived from radial basis functions (RBF-FDs) has been proposed \cite{petras2018explicit}.
This method, namely RBF-CPM, uses a smaller computational tube than the original CPM and is particularly flexible with
respect to the choice of points used to form the finite difference stencils.
In this paper, we introduce a least-squares implicit formulation of the RBF-CPM.
In contrast to previous works on implicit closest point methods on Cartesian grids \cite{macdonald2009implicit,von2013embedded},
we stabilize the method by enforcing the constant-along-normal extension using an extra equation,
and solve the resulting system by a least-squares approach.
We test our method, the {\it least-squares implicit RBF-CPM}, on a variety of examples to illustrate its convergence properties.
Examples include cases where the computational tube surrounding the surface contains inactive (and unused)
grid points.

In a second focus of this paper,
we couple the proposed method with the grid based particle method (GBPM) \cite{leung2009grid}
and perform numerical experiments on PDEs on moving surfaces. 
An extensive literature review reveals that the RBF-FD method has not been used for the solution of PDEs on moving surfaces before.
Because the least-squares implicit RBF-CPM is flexible with respect to stencil choice and exhibits good numerical stability,
its use leads to a coupled method that is robust and easily implemented.
This is an improvement over a previous combination of the original closest point method and the GBPM \cite{petras2016pdes} that
resorted to the introduction of a surface reconstruction algorithm to fill in values at deactivated nodes.

The paper unfolds as follows. In Section~\ref{ProblemSection}, we state the surface PDE under consideration
and we review the closest point method and the grid based particle method.
Section~\ref{ImplicitCPMSection} introduces the least-squares implicit closest point method using RBF-FDs and
Section~\ref{NumericalExamplesSection} presents numerical examples for static and moving surfaces.
Finally Section~\ref{SummarySection} summarizes the results and states some future work directions.

\section{Problem statement and numerical methods review}\label{ProblemSection}

In this section, we state the PDE-on-surface model under consideration and we briefly review
the closest point method and the grid based particle method (GBPM).

\subsection{Notation and formulation of PDE}\label{PDEsSection}

Following \cite{dziuk2013finite}, the conservation law of a scalar quantity $u$ with a diffusive flux on a moving surface $\Gamma(t)$ has the form
\begin{equation}\label{AdvDifPDE}
  u_t+\mathbf{v}\bigdot\nabla u+u\nabla_\Gamma\bigdot\mathbf{v}-\nabla_\Gamma\bigdot \mathcal{D}\nabla_\Gamma u = f\qquad \textrm{on } \Gamma(t),
\end{equation}
with $\mathbf{v}$ being the velocity of the surface and $\mathcal{D}$ a positive constant. If $\mathbf{n}$ is the unit normal vector of the surface $\Gamma(t)$ at some time $t$, then the velocity can be split into normal and tangential components
$\mathbf{v}=V\mathbf{n}+\mathbf{T}$, where $V = \mathbf{v}\bigdot\mathbf{n}$ and $\mathbf{T} = \mathbf{v} - V\mathbf{n}$. Using this formulation of the velocity, Equation (\ref{AdvDifPDE}) takes the form
\begin{equation}\label{AdvDifPDE2}
  u_t+V\frac{\partial u}{\partial n}-V\kappa u+\nabla_\Gamma\bigdot(u\mathbf{T}) -\nabla_\Gamma\bigdot \mathcal{D}\nabla_\Gamma u = f\qquad \textrm{on } \Gamma(t),
\end{equation}
where $\kappa$ is the mean curvature of the surface.

For a flux of Cahn-Hilliard type, the surface diffusion term in equations~(\ref{AdvDifPDE}) and (\ref{AdvDifPDE2}) is replaced by
$$\frac{1}{P_e}\nabla_\Gamma\bigdot \nu(u)\nabla_\Gamma\left(-C_n^2\Delta_\Gamma u+\frac{\partial g}{\partial u}(u)\right),$$
where $\nu$ is the mobility, $C_n$ is the Cahn number, $P_e$ is the surface Peclet number and $g$ is a double-well potential function. Consequently, we find
\begin{equation}\label{CahnHilPDE}
  u_t+V\frac{\partial u}{\partial n}-V\kappa u+\nabla_\Gamma\bigdot(u\mathbf{T}) - \frac{1}{P_e}\nabla_\Gamma\bigdot \nu(u)\nabla_\Gamma\left(-C_n^2\Delta_\Gamma u+\frac{\partial g}{\partial u}(u)\right)
= f\qquad \textrm{on } \Gamma(t).
\end{equation}

\subsection{The closest point method}\label{CPMsubsection}
The closest point method \cite{ruuth2008simple} is an embedding method for solving PDEs on static smooth surfaces.
Given a uniform Cartesian grid that contains a surface $\Gamma$, the closest point function maps each grid point to its closest point on the surface:

\begin{definition}\label{ClosestPointRepresentation}
  Let $\mathbf{z}$ be some point in the embedding space $\mathds{R}^d$ that is sufficiently close to $\Gamma$. Then,
  $$cp_\Gamma(\mathbf{z})=\arg\min_{\mathbf{x}\in \Gamma}\|\mathbf{x}-\mathbf{z}\|_2$$
  is the closest point to $\mathbf{z}$ on the surface $\Gamma$.
\end{definition}

To compute the closest point function, we select a method appropriate for the surface under consideration. For simple surfaces such as the circle, sphere, and torus, we typically use an analytical formula, e.g.,  $cp_\Gamma(\mathbf{z})= r \cdot  \mathbf{z}/\|\mathbf{z}\|_2$ for a sphere of radius $r$ centered at the origin. For parametrized surfaces (e.g., an ellipse, ellipsoid or M\"obius strip), we compute the closest point function by minimizing distance over the free parameters; see, e.g., \cite{merriman2007diffusion} for results based on this technique. The other surface representation that we encounter frequently is triangulated form. Here, we compute the closest point function by looping over the list of triangles according to the algorithm provided in \cite{macdonald2008level}. In this approach, for each grid node in a suitable neighborhood of a triangle $T_i$, the closest point on $T_i$ is computed and stored. After looping through all the triangles, the closest point on the surface for any grid node is simply the closest point over all stored possibilities. See \cite{macdonald2008level} for details on this procedure.

The grid nodes and the corresponding closest point values together form a {\it closest point representation}. Note that the normal to the surface is not required or computed in the classical closest point method.

The closest point representation maps each grid point to its closest
point on the surface. By replacing grid node values by values at the closest point, we obtain a constant normal extension
of the surface values.
The surface PDE is extended into the embedding space by replacing derivatives intrinsic to the surface with the corresponding Cartesian derivatives according to two principles \cite{ruuth2008simple}:

\begin{principle}\label{equivalence of gradients}
  Let $v$ be a function on $\mathds{R}^d$ that is constant along normal directions of $\Gamma$.
  Then, at the surface, intrinsic gradients are equivalent to standard gradients, $\nabla_\Gamma v=\nabla v$.
\end{principle}

\begin{principle}\label{equivalence of divergence}
  Let $\mathbf{v}$ be a vector field on $\mathds{R}^d$ that is tangent to $\Gamma$ and tangent to all surfaces displaced
  by a fixed distance from $\Gamma$. Then, at the surface, $\nabla_\Gamma\bigdot\mathbf{v}=\nabla\bigdot\mathbf{v}$.
\end{principle}
By combining Principles~\ref{equivalence of gradients} and \ref{equivalence of divergence}, higher-order surface derivatives, such as the Laplace-Beltrami and biharmonic operators,
can be replaced with the corresponding standard Cartesian derivatives in the embedding space \cite{ruuth2008simple,macdonald2009implicit}.  
For a variety of related theory, see \cite{marz2012calculus, cheung2018kernel}.

The algorithm of the closest point method alternates two steps: the extension of the solution into the embedding space,
and the solution of the PDE.
The extension is an interpolation step, as the points on the surface are not necessarily grid points. Taking into consideration the size of the
interpolation and differencing stencils, a computational tube around the surface can be formed.
For a second-order finite difference approximation of the Laplace-Beltrami operator,
we may choose the computational tube radius to be
\begin{equation}\label{Bandwidth}
  \gamma_{CPM}=\sqrt{(d-1)\left(\frac{p+1}{2}\right)^2+\left(1+\frac{p+1}{2}\right)^2}\Delta x
\end{equation}
in the $d$-dimensional embedding space uniformly discretized using a Cartesian grid with spatial step-size $\Delta x$, where $p$ is the degree of the interpolating polynomial \cite{ruuth2008simple}.

An implicit closest point method that allows the use of large time step-sizes in implicit time discretization schemes is also available \cite{macdonald2009implicit}.
To illustrate, assume that $\Delta_\Gamma$ is the Laplace-Beltrami operator on a surface $\Gamma$ and that $cp_\Gamma$
gives its closest point representation. Then, following Principles \ref{equivalence of gradients} and
\ref{equivalence of divergence}, $\Delta_\Gamma$ is replaced with the Cartesian differential operator $\Delta$ in the embedding PDE, i.e.,
\begin{equation}\label{Surface operator}
  \Delta_\Gamma u = \Delta u(cp)
\end{equation}
on the surface, where $u$ is a scalar function. In discretized form, for a uniform Cartesian grid that contains the surface $\Gamma$, Equation~(\ref{Surface operator}) is approximated as
$$\Delta_hEU=:\widetilde{M}U$$
where $\Delta_h$ is the discretized differential operator using finite differences (e.g., second-order centered finite differences), $E$
is the closest point extension matrix (typically formed via barycentric Lagrange interpolation) and $U$ is the discretized solution $u(cp)$. Unfortunately, the matrix $\widetilde{M}$
has eigenvalues with positive real components, leading to instability \cite{macdonald2009implicit}.
Stable computations are obtained by discretizing $\Delta$ using
$$M = diag(\Delta_h)+(\Delta_h-diag(\Delta_h))E.$$
Note that the matrix $M$ approximates \eqref{Surface operator} as $\widetilde{M}$ on the surface, but has all its eigenvalues
in the left half of the complex plane.

\subsection{The grid based particle method}\label{GBPMSection}
To evolve the surface based on a velocity $\mathbf{v}$, we use the grid based particle method (GBPM) by Leung and Zhao \cite{leung2009grid}.
To initialize the GBPM, a Cartesian grid is constructed that contains the surface $\Gamma$. Over a neighborhood of the surface of radius $\gamma_{GBPM}$, called the {\it computational tube}, a closest point representation of $\Gamma$ is constructed. Specifically, the grid points contained in the computational tube are mapped to their closest points on the surface. In the GBPM, this mapping might not be constructed for all the grid points within the neighborhood of the surface \cite{petras2016pdes}. The grid points
 that are contained within the computational tube and are mapped to their closest points are called \emph{active grid points} and their closest points on the
 surface are called \emph{footpoints}.

Following the initialization described in the previous paragraph, the system is evolved in time.   Each time step of size $\Delta t$ consists of three steps:
\begin{figure}
    \centering
    \includegraphics[width=0.32\textwidth]{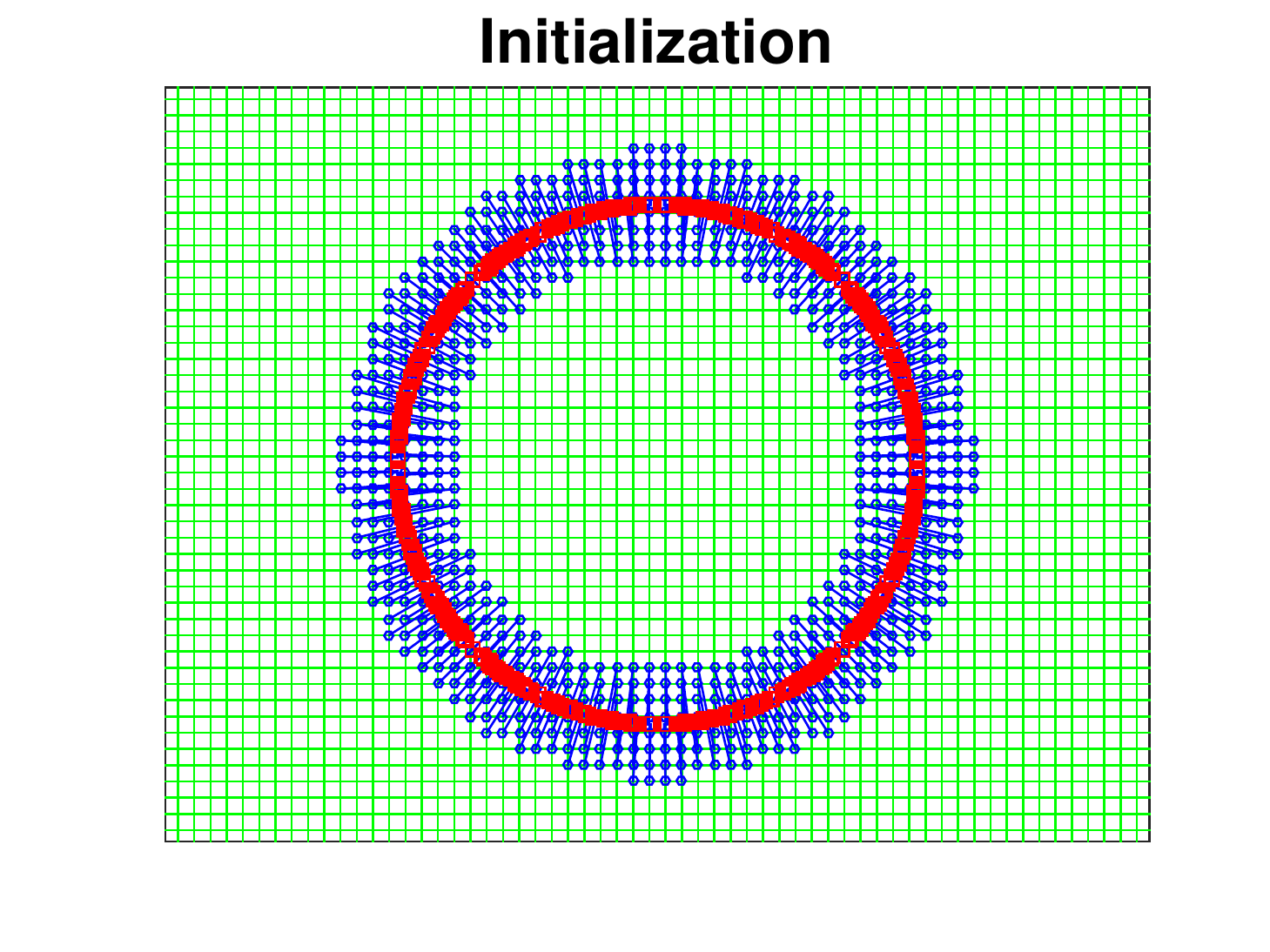}
    \includegraphics[width=0.32\textwidth]{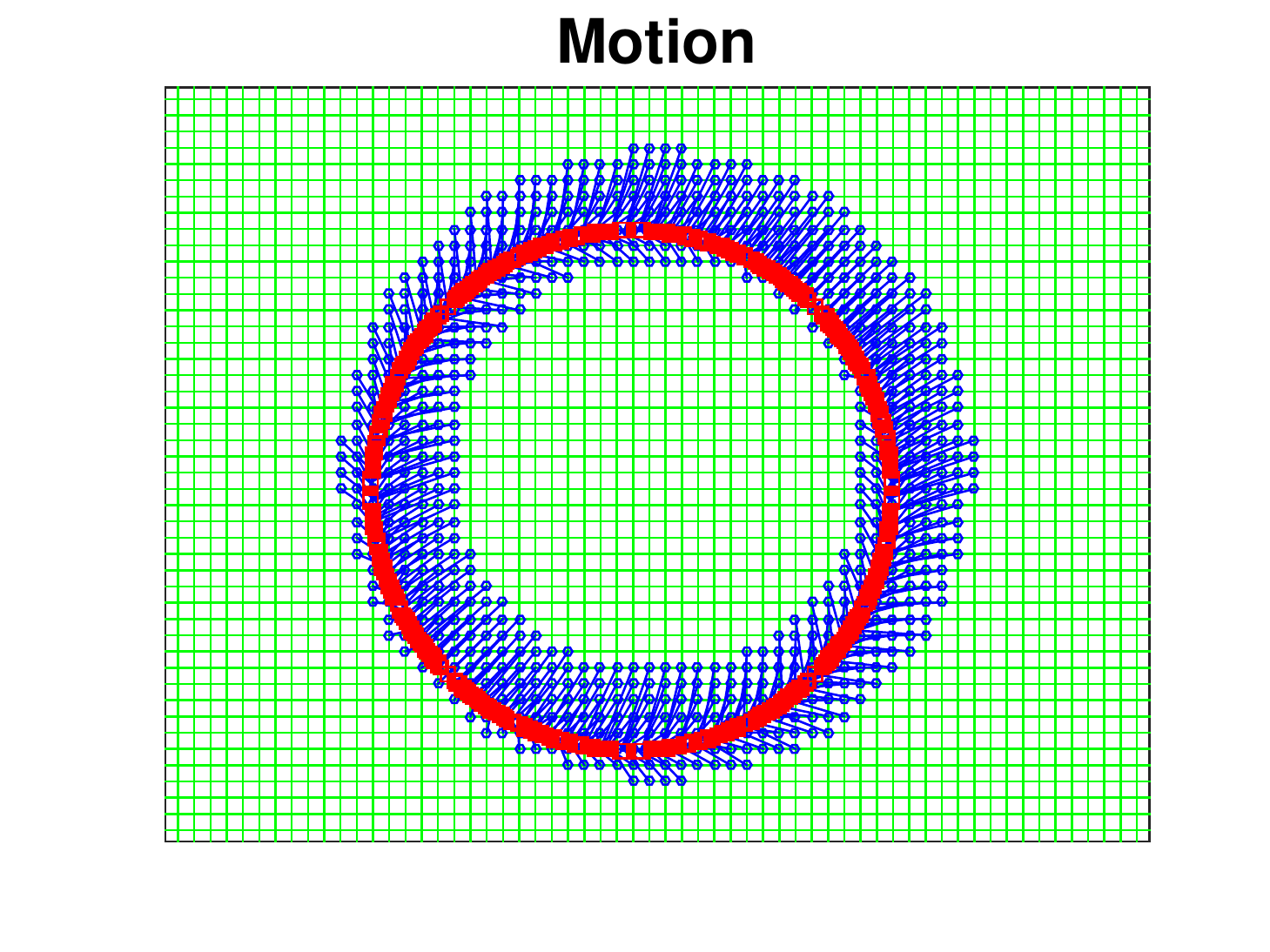}
    \includegraphics[width=0.32\textwidth]{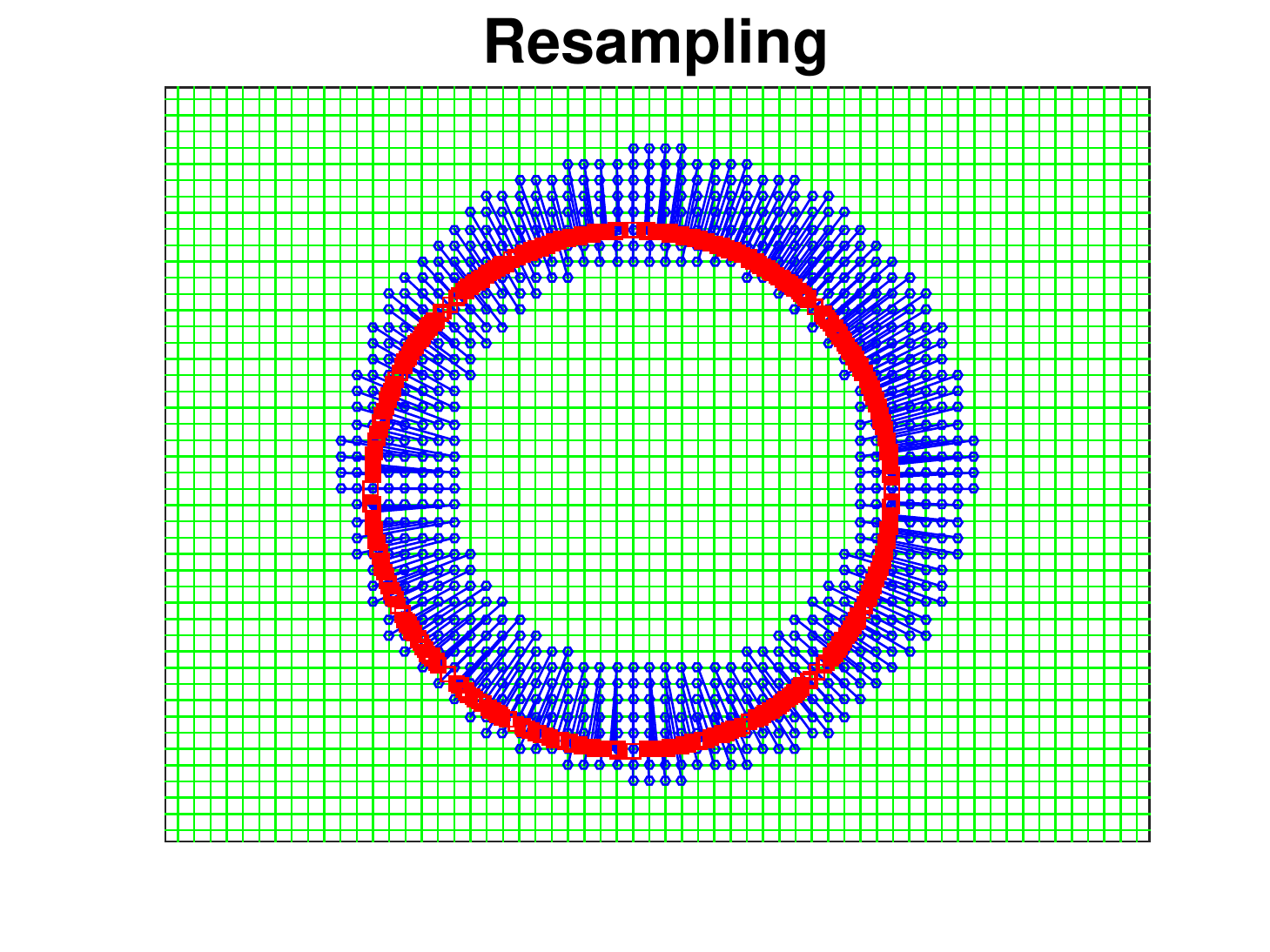}\\
    \includegraphics[width=0.33\textwidth]{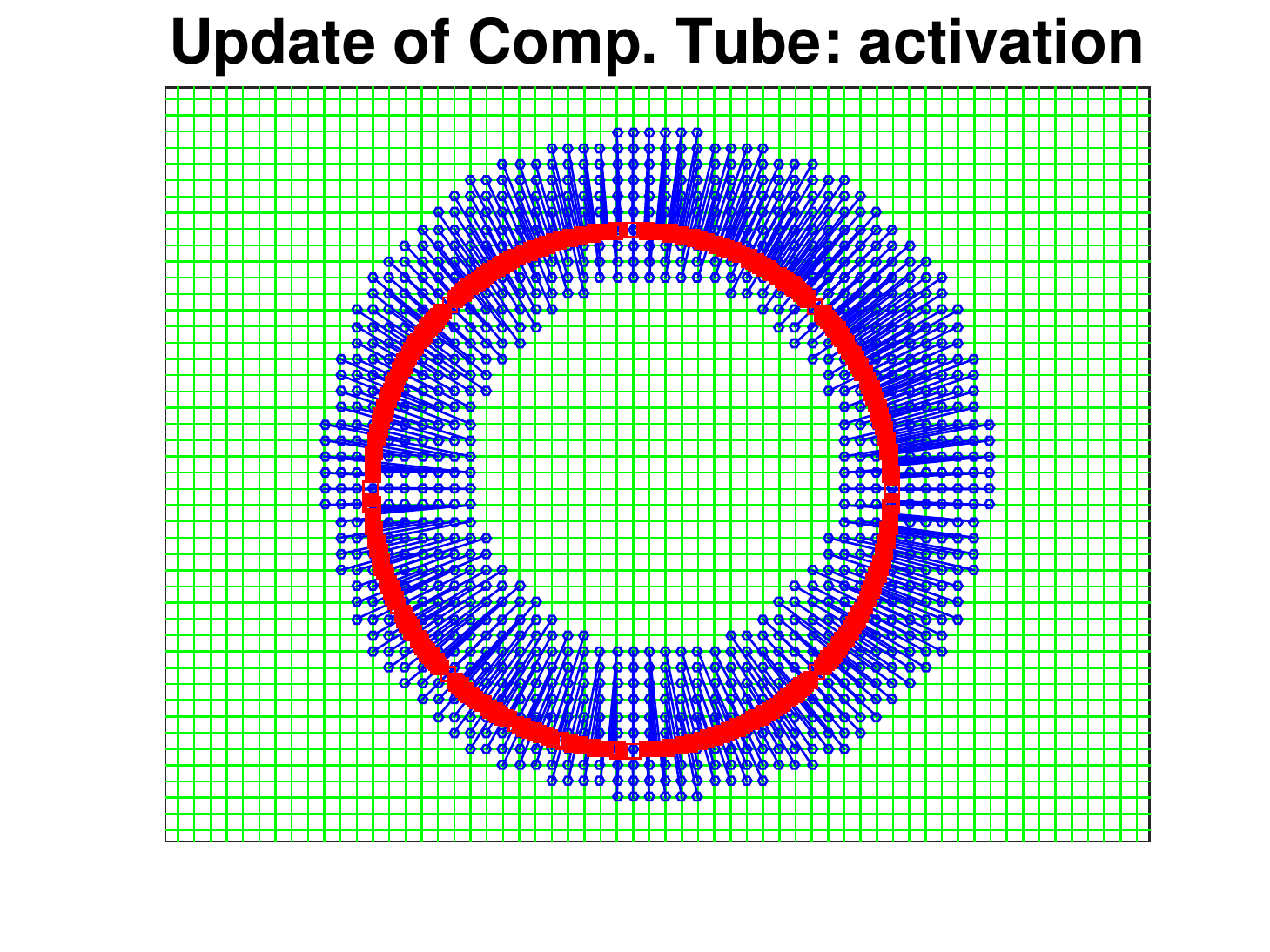}
    \includegraphics[width=0.33\textwidth]{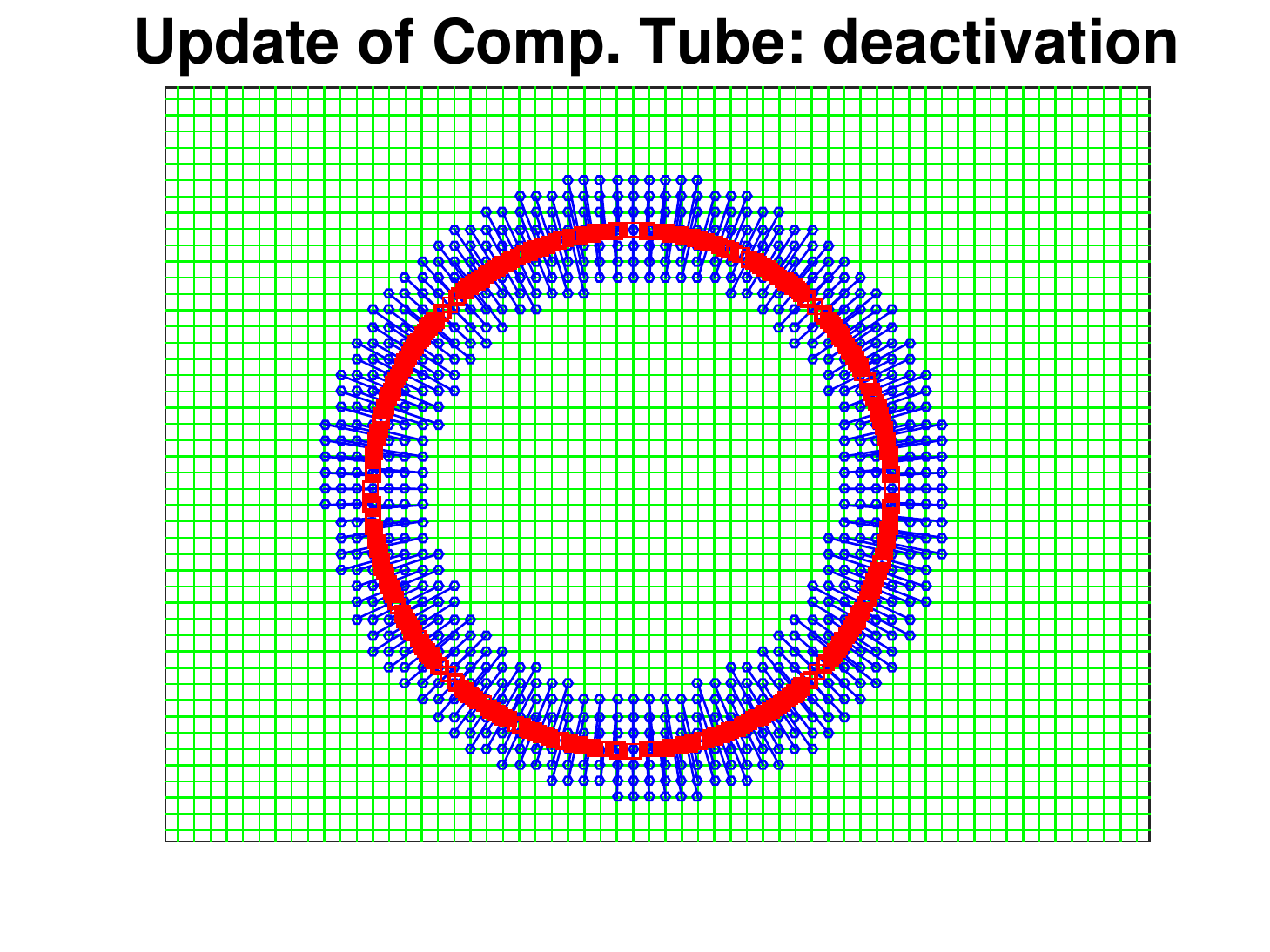}
    \caption{An illustration of the main steps of the GBPM (from top left to bottom right). Active grid points (blue dots)
    are connected to their footpoints (red dots) with blue lines. The green lines correspond to the grid.}
    \label{GBPM}
\end{figure}

\begin{enumerate}
   \item \textbf{Motion}: The footpoints are moved according to a motion law.
  \item \textbf{Resampling}: For each active grid point, a new closest point to the surface (as defined by the footpoints) is computed.  This gives the updated footpoints.
  \item \textbf{Update of the Computational Tube}: There are two stages in this step. During the first stage, all the grid points
  that have neighboring active grid points are activated and the resampling step is applied to find their footpoints.
  During the second stage, all the grid points whose distance from their footpoints is larger than the tube radius $\gamma_{GBPM}$ are deactivated.
\end{enumerate}
Figure~\ref{GBPM} illustrates the main steps of the GBPM algorithm. We shall use the GBPM to evolve closed surfaces by curvature-dependent motions, however the method is also capable of capturing the motion of open surfaces \cite{leung2009gridopen}
and higher-order geometric motions of surfaces \cite{leung2011grid}.

\section{A least-squares implicit RBF closest point method}\label{ImplicitCPMSection}
In this section, we introduce a least-squares implicit RBF closest point method (RBF-CPM).  Numerical experiments
are provided to illustrate the convergence and behavior of the method for a variety of problems
on static surfaces.

\subsection{Description of the method}\label{ImplicitCPMDescription}
In \cite{petras2018explicit}, an RBF-CPM method is presented.
Numerical results show the method's potential. A notable advantage of the method is the reduction of
the size of the computational tube (e.g., by $25\%$ relative to the standard finite difference implementation
of the closest point method for the case of the Laplace-Beltrami operator).
The method also gains flexibility with respect to the stencil choice due to its use of RBF-FD; in particular, high-order accuracy can be achieved simply by increasing the number of points in the RBF-FD stencil.
On the other hand, the method can be very slow when applied to stiff problems due to its use of explicit time stepping.
We therefore seek an implicit formulation of the RBF-CPM.

For illustration purposes, consider the heat equation intrinsic to a surface $\Gamma$,
\begin{equation}\label{HeatEquation}
  u_t=\Delta_\Gamma u.
\end{equation}
Let $\Omega$ be the embedding space of the surface $\Gamma$. We denote the constant-along-normal extension by $\tilde{u} = u(cp_\Gamma):\Omega\rightarrow\mathbb{R}$ so that
\[\tilde{u}(\mathbf{x}) = u(\mathbf{x}),\text{ for all }\mathbf{x}\in\Gamma.\]
Equation~(\ref{HeatEquation}) can be written in an equivalent form in the embedding space $\Omega$ as
\begin{equation}\label{HeatEquationExtended}
  \tilde{u}_t=\Delta \tilde{u}.
\end{equation}
To discretize \eqref{HeatEquationExtended}, a Cartesian grid is constructed in a tubular neighborhood around the surface in the embedding space (i.e., the grid shown in blue in Figure~\ref{GBPM}), $Z = \{\mathbf{z}_j\}_{j=1}^{n_Z}$. Using a closest point representation of the surface $cp_\Gamma$ as in Definition~\ref{ClosestPointRepresentation}, we define the surface points $X=\{\mathbf{x}_j\}_{j=1}^{n_Z}\subset \Gamma$, such that $cp_\Gamma(Z) = X$.

Following \cite{petras2018explicit}, for a surface point $\mathbf{x}_j = cp_\Gamma(\mathbf{z}_j)$, $\mathbf{z}_j\in\Omega$, and a collection of $m$ grid points $Z_j=\{\mathbf{z}_{j_1},...,\mathbf{z}_{j_m}\}\subset Z$ closest to $\mathbf{x}_j$, a local RBF-FD approximation gives
\begin{equation}
\Delta\tilde{u}(\mathbf{x}_j) \approx \mathbf{w}_j\tilde{u}(Z_j),
\label{eq:RBF-FD weights}
\end{equation}
where $\mathbf{w}_j\in \mathbb{R}^{1\times m}$.

We consider the use of polyharmonic spline (PHS) RBFs $\phi(r) = r^{2k-1}$, for integer $k > 1$, for the construction of the RBF-FD approximations, however other RBFs can be used. The local RBF-FD weights are calculated as
$$\tilde{\mathbf{w}}_j = (B(\mathbf{x}_j,Z_j)A(Z_j,Z_j)^{-1})^T,$$
and $\mathbf{w}_j = (\tilde{w}_{j_1},...,\tilde{w}_{j_m})$ consists of the first $m$ terms of $\tilde{\mathbf{w}}_j$.
The matrix $A$ and the vector $B$ are expressed as in \cite{petras2018explicit}, augmented with polynomial terms as described in \cite{flyer2016role}. Specifically, the matrix $A$ and the vector $B$ have the form

\begin{gather*}
\resizebox{\textwidth}{!}{$A(Z_j,Z_j) = \left(\begin{array}{ccccccc}
 &  &  & \vline &  &  & \\
 & \phi(\|\mathbf{z}_{j_i} - \mathbf{z}_{j_l}\|) &  & \vline &  & Q & \\
 &  &  & \vline &  &  & \\
\hline 
 & Q^T & & \vline &  & 0\\
\end{array}\right),
B(\mathbf{x}_j,Z_j) = \left(\begin{array}{c}
\Delta\phi(||\mathbf{x}_j-\mathbf{z}_{j_1}||)\\
\vdots\\
\Delta\phi(||\mathbf{x}_j-\mathbf{z}_{j_m}||)\\
\hline 
 R\\
\end{array}\right),$}
\label{eq:mateqRBFCPM}
\end{gather*}
where $Q :=[q_i(\mathbf{z}_{j_l})]$, $R = [(\Delta q_i)(\mathbf{x}_j)]$ and $q_i$ are basis of the polynomial space $\mathbb{P}_p$.

From the local systems described above, a global sparse matrix can be introduced as in  \cite{petras2018explicit}, i.e.

\begin{equation}\label{Wu}
  \Delta_\Gamma u(X) =  \Delta \tilde{u}(X) \approx W \tilde{u}(Z),
\end{equation}
and (\ref{HeatEquationExtended}) yields the ODE system
\begin{equation}\label{SemiDiscreteHeatEquation}
  \dot{U}_Z=W U_Z,
\end{equation}
where $U_Z$ is the semi-discretized solution in the embedding space.

To obtain an improvement in efficiency over explicit methods,
we applied backward Euler and other implicit time stepping methods
to (\ref{SemiDiscreteHeatEquation}).
Unfortunately, in numerical experiments, these schemes did not yield an improvement in the observed stability time step restriction. In particular, using backward Euler for the time discretization of equation~\eqref{SemiDiscreteHeatEquation} provided a stable solution for small time step-sizes, similar to forward Euler, yet instabilities were observed for large time step-sizes.
To stabilize~(\ref{SemiDiscreteHeatEquation}), one can consider enforcing the constant normal extension of
the solution as part of the ODE system (cf. \cite{macdonald2009implicit,Macdonald20117944,von2013embedded}),
$$\dot{U}_Z =WU_Z-c(U_Z-PU_Z),$$
where $c$ is a positive constant and $P$ is the projection matrix from the embedding space $\Omega$ on the surface $\Gamma$ \cite{petras2018explicit}, i.e.

\begin{equation}\label{Pu}
 u(X) = \tilde{u}(X) \approx P \tilde{u}(Z),
\end{equation}
which can be constructed similarly to the matrix $W$, by replacing the Laplacian with the identity map in (\ref{eq:RBF-FD weights}). The matrix $P$ is the corresponding interpolation matrix or extension matrix $E$ in the classical closest point method \cite{ruuth2008simple}, and is constructed by local RBF-FD systems. Yet, the identification process of the proper constant $c$ that increases the stability time step restriction over implicit schemes is unclear.  

Inspired by \cite{piret2012orthogonal,von2013embedded}, we propose an alternative approach for stabilizing the method.
Specifically, we enforce the solution $\tilde{u}$ of \eqref{HeatEquationExtended} to be a constant-along-normal extension of $u$ by introducing an extra equation,
$\tilde{u} = u(cp_\Gamma)$. This extra equation should hold for all times $t$. The semi-discretized system of equations becomes
\begin{equation}\label{SemiDiscreteHeatEquationLeastSquares}
  \left\{
    \begin{array}{l}
      \dot{U}_Z=W U_Z,\\
      U_Z - PU_Z = 0. \\
    \end{array}
  \right.
\end{equation}
Using the method of lines approach and 
applying an implicit time stepping method in \eqref{SemiDiscreteHeatEquationLeastSquares} leads to an over-determined system.
For example, using the techniques described in \cite{petras2018explicit}, the application of the backward differentiation formula BDF2 yields
\begin{equation}\label{DiscreteHeatEquationLeastSquares}
  \left\{
    \begin{array}{l}
      U^{n+1}_Z = \left(\frac{4}{3}U^n_Z - \frac{1}{3}U^{n-1}_Z\right)+\frac{2}{3}\Delta tW U^{n+1}_Z,\\
      U^{n+1}_Z - PU^{n+1}_Z = 0, \\
    \end{array}
  \right.
\end{equation}
where $U^n_Z$ is the discretized approximate solution at time $n \Delta t$. For notational simplicity, let us introduce a matrix $A$ that contains all the coefficients of the implicit terms, i.e., $A = I-2/3\Delta tW$, where $I$ is the identity matrix. Further, denote by the vector $b^n$ all the explicit terms, i.e., $b^n = 4/3U^n_Z - 1/3U^{n-1}_Z$. Thus, the equations above take the form $AU^{n+1}_Z = b^n$, $(I-P)U^{n+1}_Z=0$. There is no solution, in general, that satisfies both equalities; hence we consider the solution of the minimization problem
\begin{equation}\label{eq:discreteLS}
U^{n+1}_Z = \arginf_{w\in\mathbb{R}^{n_Z}}\left(\|Aw-b\|^2_{\ell_2} + c\|(w-Pw\|^2_{\ell_2}\right),
\end{equation} 
where $c>0$ is a penalty constant. Instead of $\ell_2(\mathbb{R}^{n_Z})$, one can use other norms in \eqref{eq:discreteLS}, or even mixed norms. However, since we are dealing with smooth solutions of some diffusion equations, any $\ell_p$ norm will behave similarly. An advantage of selecting $\ell_2$ is that it approaches the $L_2(\Omega)$ norm up to some scaling factors, as the size of the uniform grid $h_Z\rightarrow0$.

Using the least squares method, the solution of the minimization problem \eqref{eq:discreteLS} takes the form
$$ U^{n+1}_Z = \left(\begin{array}{c}
A\\
c(I-P)
\end{array}\right)^+\left(\begin{array}{c}
b\\
0
\end{array}\right),$$
where $+$ denotes the left pseudoinverse of the matrix.
In this work, we fix $c=1$ because we have found that the solution quality is robust with respect to changes in $c$. This approach was found to be unconditionally stable in practice. See \cite{cheung2018kernel} for some background on the convergence and the oversampling requirements of the least-squares method with RBFs for elliptic equations.


The use of the PHS RBF in the calculation of the local RBF-FD stencils described in \cite{petras2018explicit} as well as the least-squares method provide flexibility in the node placement and the computational tube regularity.
Of particular interest is the case where the closest point mapping to the surface is not available for all the grid points in a neighborhood of the surface, thus introducing irregular computational tube patterns \cite{petras2016pdes}. 
We shall see (in Section~\ref{sec:numex} below) that the proposed stabilization is particularly effective
for computational tubes of this type.

\subsection{Numerical experiments} \label{sec:numex}
In this section, we present numerical experiments for the solution of PDEs on static surfaces. For the solution
of the sparse matrix least-squares systems, we use the MATLAB code Factorize\footnote{Available at
http://www.mathworks.com/matlabcentral/fileexchange/24119-don-t-let-that-inv-go-past-your-eyes--to-solve-that-system--factorize-} \cite{davis2013algorithm}.
Unless stated otherwise, the PHS RBF $\phi(r) = r^7$ is used with augmented polynomial basis that span $\mathbb{P}_3$. The number of points in the RBF-FD stencil $m$ is chosen as twice the number of the augmented polynomial terms, as suggested in \cite{flyer2016role}. The computational tube radius $\gamma$ is chosen according to the Gauss circle problem, as described in \cite{petras2018explicit}. 

For computational tubes with deactivated/missing grid nodes, changes in the tube radius $\gamma$ are unnecessary. Specifically, for RBF-FD stencils that consist of the closest grid points, a suitable stencil can be found using $m$ closest neighbors to a surface point. The least-squares approach stabilizes for such irregular stencils.

\subsubsection{Heat equation on a circle}
\begin{figure}[h!]
\centering
\includegraphics[width = 0.5\textwidth]{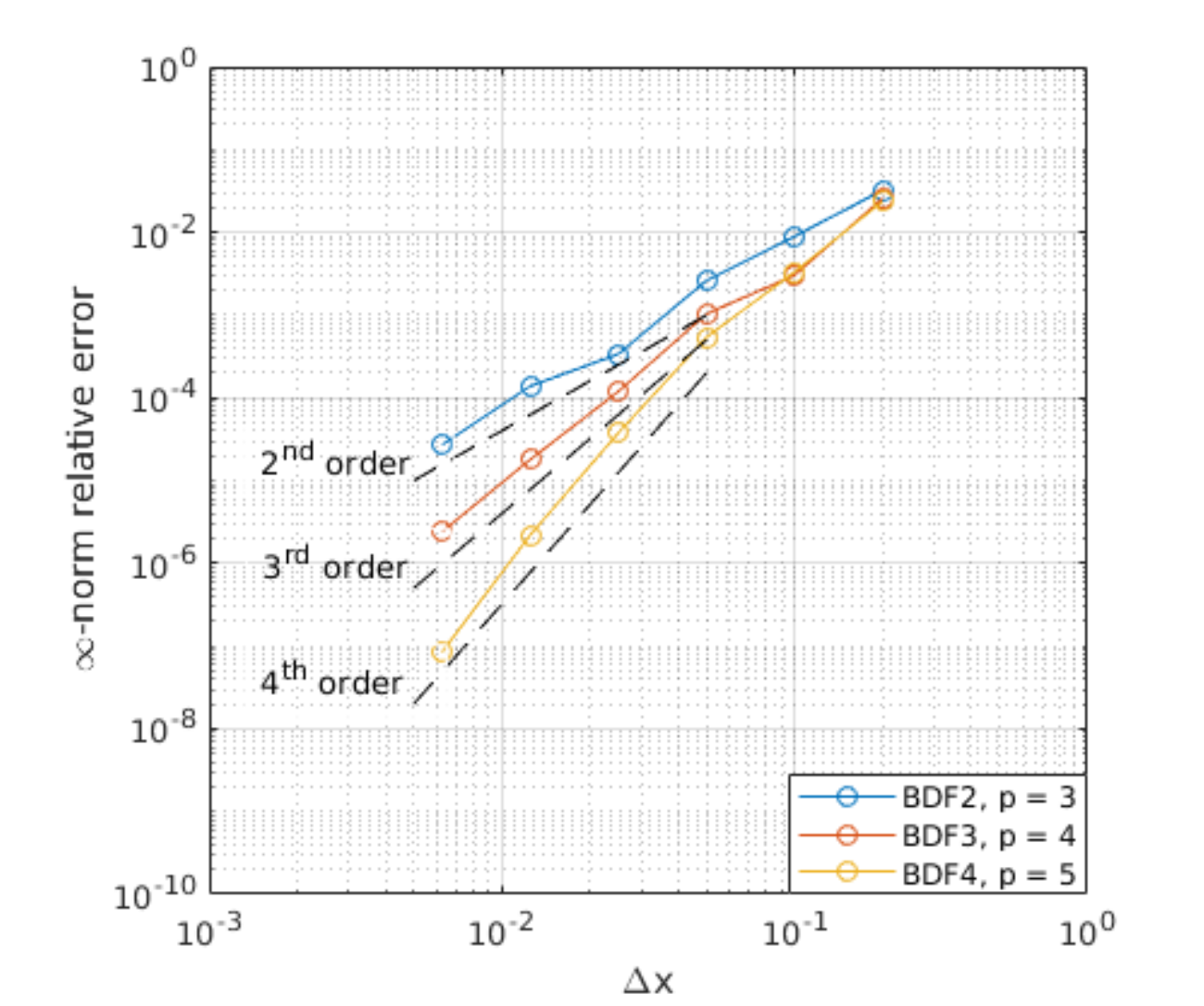}
\caption{The relative error of the heat equation on the unit circle at $t=1$ using different backward differentiation formulas and augmented polynomial degree in our RBF-FD calculations.}\label{figure: BDF convergence}
\end{figure}

Our first example approximates the heat equation
$$u_t=\Delta_\Gamma u$$
on the unit circle using $u(\theta,0)=\sin\theta$ as the initial condition. The exact solution at all times $t>0$ is
$$u(\theta,t) = e^{-t}\sin\theta.$$
We begin by applying the backward differentiation formulas of second-order (BDF2), third-order (BDF3), and fourth-order (BDF4) to explore the convergence of the proposed method. Similar to \eqref{DiscreteHeatEquationLeastSquares}, the projection matrix $P$ is applied to all the explicit terms in the BDF discretizations. Using a time step-size of $\Delta t = \Delta x$ and the exact solution as initial steps for the BDF schemes, Figure~\ref{figure: BDF convergence} shows the $\ell_\infty$-norm error of the approximate solution relative to the exact solution at time $t=1$. Observe that the expected RBF-FD spatial convergence for different degrees $p$ of augmented polynomials (see Section~\ref{ImplicitCPMDescription}) can be achieved using the BDF discretizations of the corresponding convergence rate.

\begin{figure}[h!]
\centering
\includegraphics[width = 0.49\textwidth]{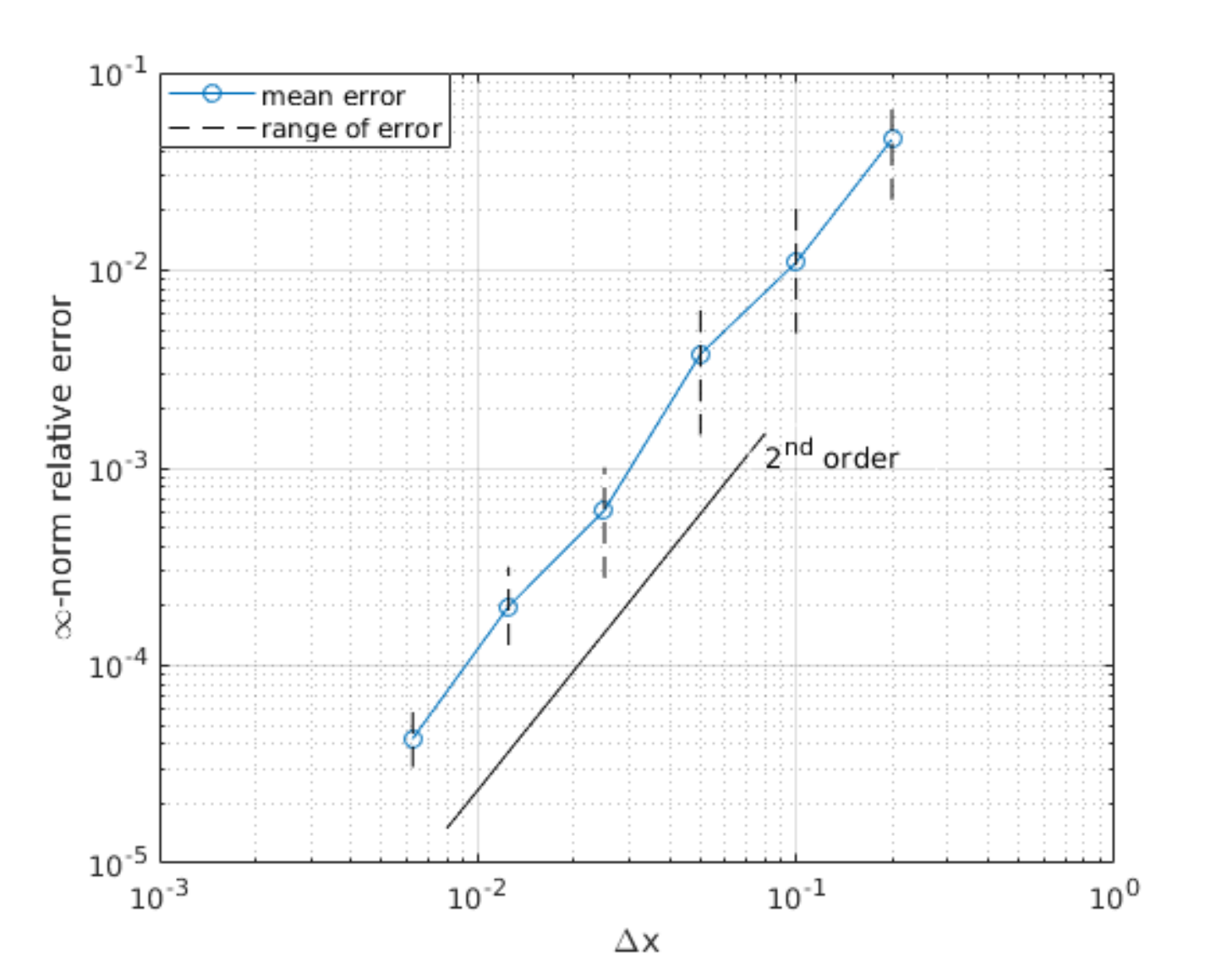}
\includegraphics[width = 0.49\textwidth]{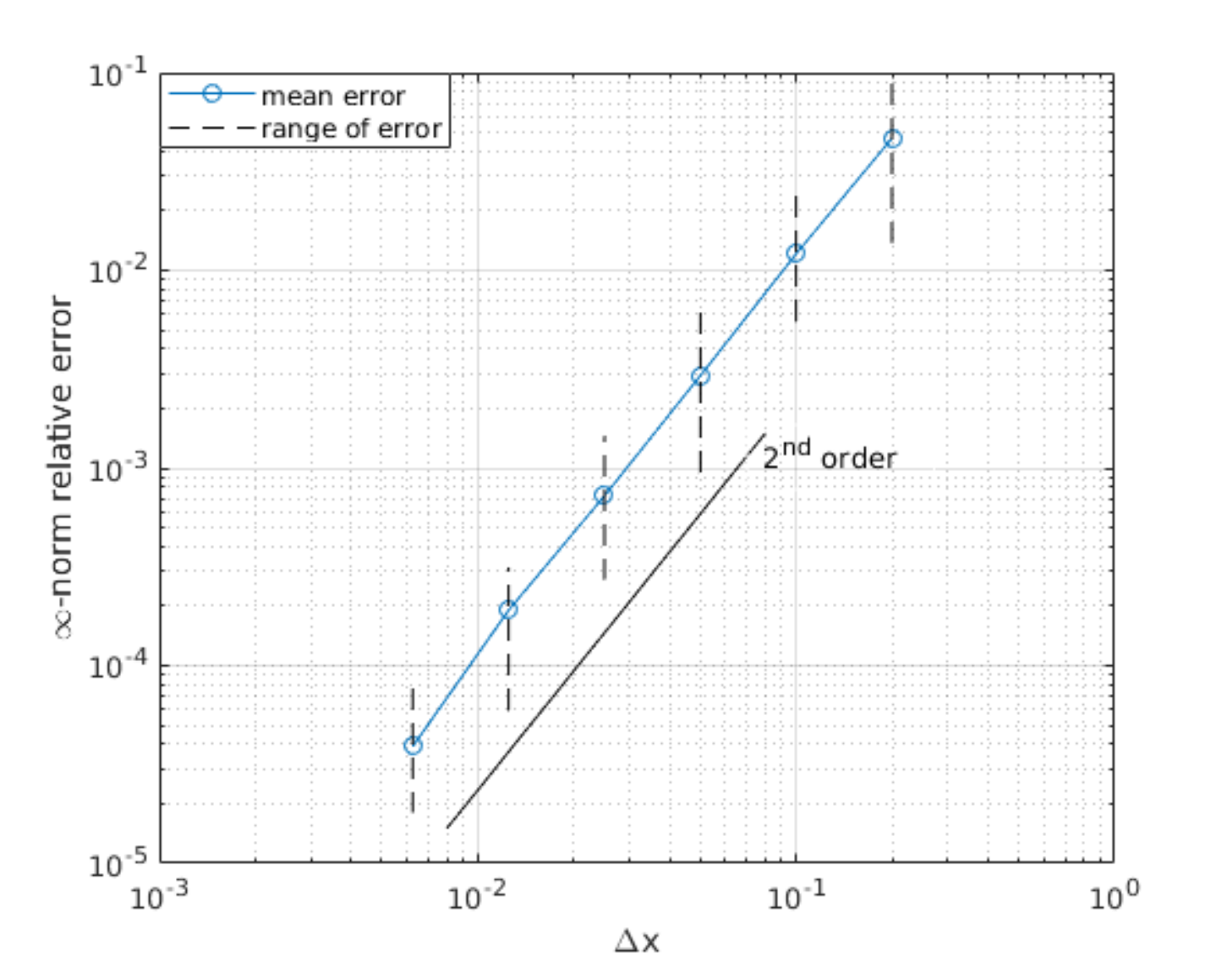}
\caption{The relative error of the heat equation on the unit circle for $1\%$ (left) and $5\%$ (right) point removal at $t=1$ using BDF2 and up to third degree augmented polynomials in the RBF-FD calculation.}\label{figure: circle holes convergence}
\end{figure}

Next, we consider the BDF2 scheme with augmented polynomial basis that spans $\mathbb{P}_3$ and explore the convergence of the method for computational tubes with holes (grid nodes that cannot be mapped to their closest points on the surface). To demonstrate the convergence of the method for computational tubes with holes, we randomly remove
$1\%$ and $5\%$ of the points in the computational tube and apply the method to the same problem. For each spatial discretization level $\Delta x$,
50 experiments are carried out by removing $1\%$ and $5\%$ of the grid points within the computational tube
(along with corresponding nodal values).
The mean value of the relative error as well as the range of the error over the 50 experiments for different grid spacings $\Delta x$ is shown in Figure~\ref{figure: circle holes convergence}.
As expected, we obtain order $p-1$ convergence by combining RBF-FDs augmented with polynomials that span $\mathbb{P}_p$, with BDF discretizations of order $p-1$.

\subsubsection{Heat equation on a sphere}
\begin{figure}[h!]
\centering
\includegraphics[width = 0.49\textwidth]{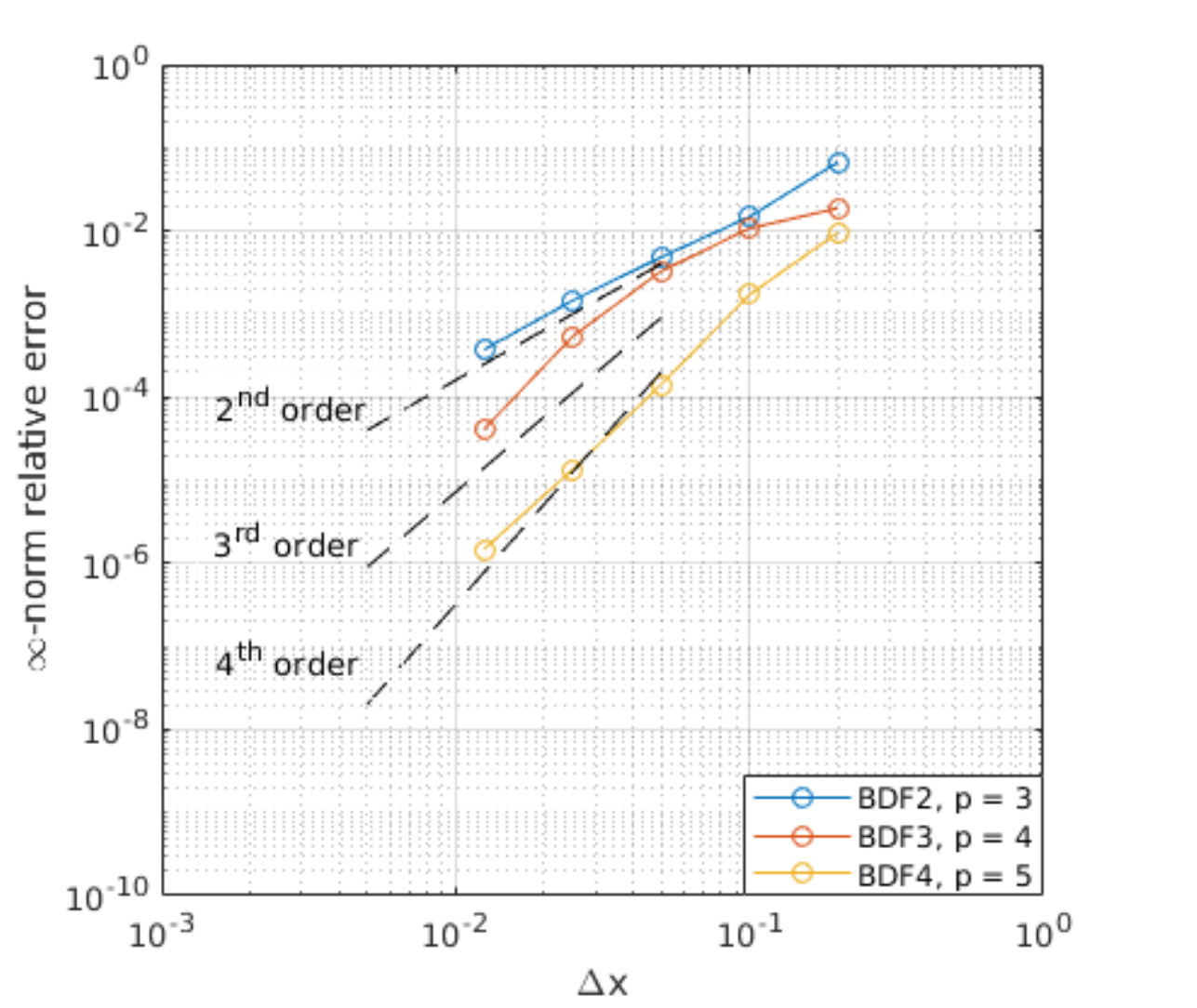}
\caption{The relative error of the heat equation on the unit circle at $t=1$ using different backward differentiation formulas and augmented polynomial degree in our RBF-FD calculations.}\label{figure: BDF convergence sphere}
\end{figure}

Consider the heat equation on the unit sphere. For the parametrization of the sphere
$$\mathbf{x}(\theta,\phi)=(\cos\theta\cos\phi,\sin\theta\cos\phi,\sin\phi)$$
and the initial profile
$$u(\theta,\phi,0)=\sin\phi,$$
the exact solution for all times $t>0$ is given by
$$u(\theta,\phi,t)=e^{-2t}\sin\phi.$$
Similarly to the example on the unit circle, we explore the convergence of the method by introducing the BDF2, BDF3 and BDF4 discretizations in time for augmented polynomial basis of $\mathbb{P}_p$, with $p=3,4$ and $5$ degree in the RBF-FD calculation. Figure~\ref{figure: BDF convergence sphere} shows the relative error at time $t=1$ for $\Delta t=\Delta x$. For RBF-FDs augmented with polynomial basis of $\mathbb{P}_p$, we see a $p-1$ approximate order of convergence in $\Delta x$, in agreement with our expectations. 

\begin{figure}[h!]
\centering
\includegraphics[width = 0.49\textwidth]{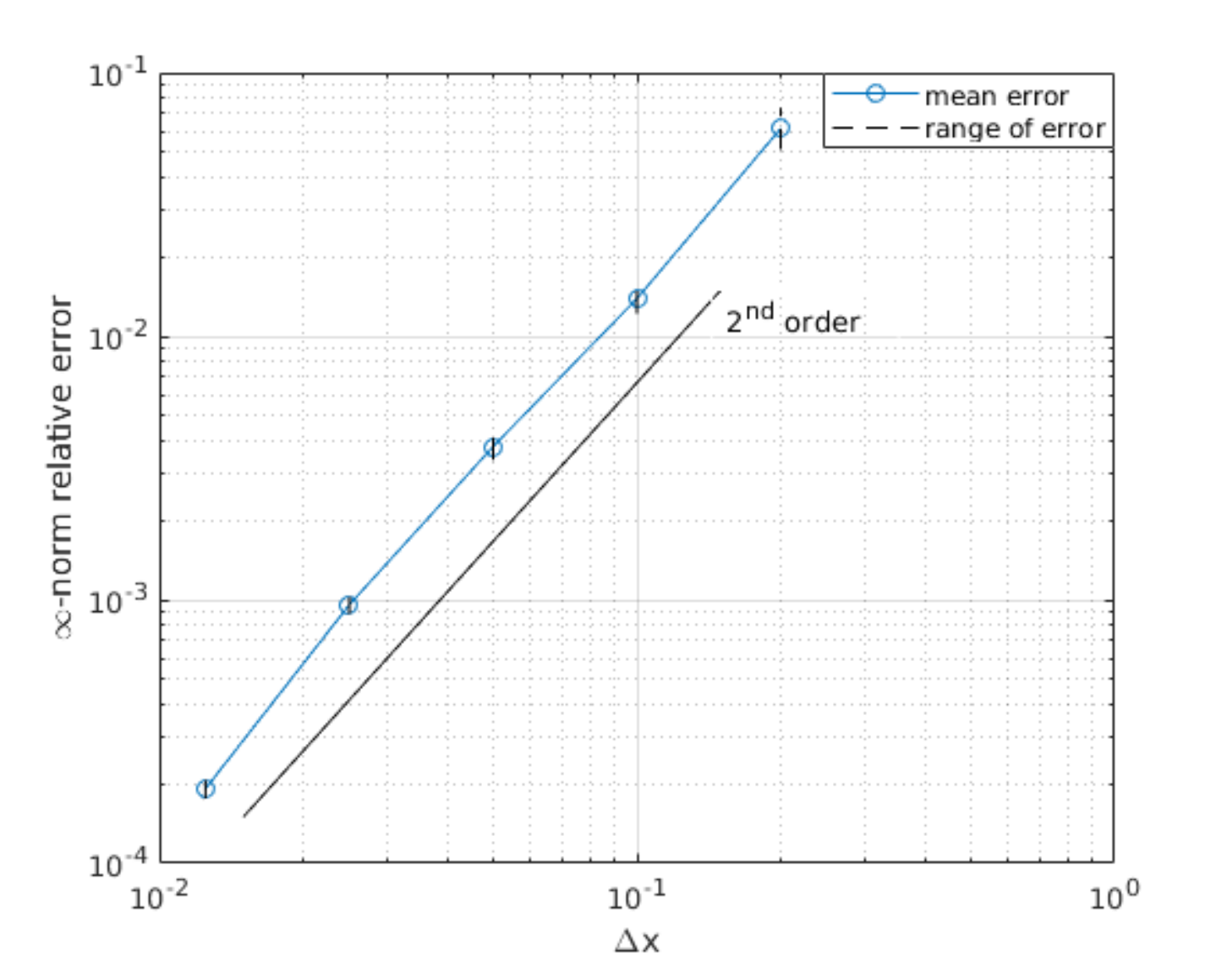}
\includegraphics[width = 0.49\textwidth]{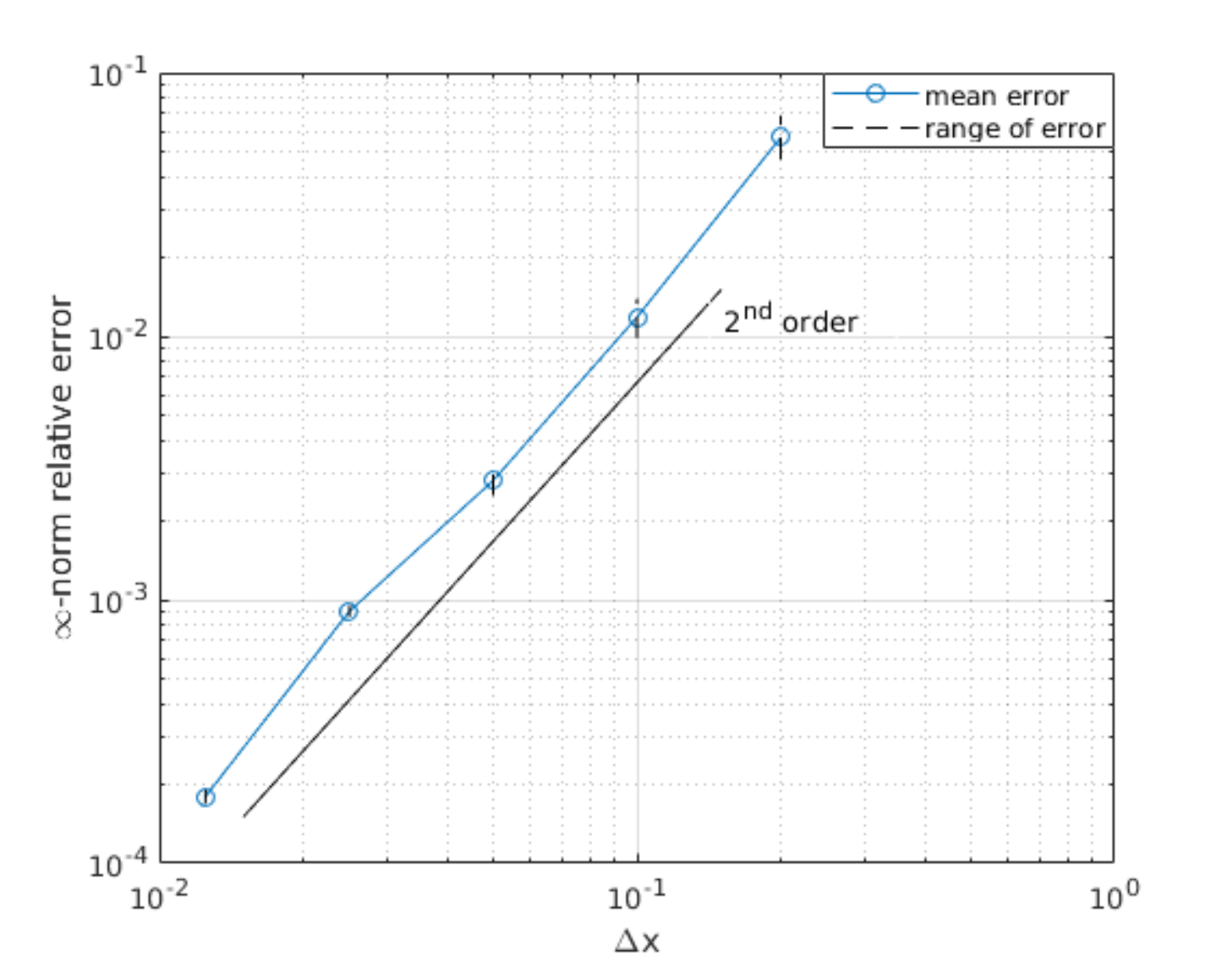}
\caption{The relative error of the heat equation on the unit sphere for $1\%$ (left) and $5\%$ (right) point removal at $t=1$ using BDF2 and up to third degree augmented polynomials in the RBF-FD calculation.}\label{figure: sphere holes convergence}
\end{figure}

We also perform a convergence analysis on computational tubes with holes by randomly removing
$1\%$ and $5\%$ of the points and approximating the solution of the same problem using BDF2 and $p=3$. Figure~\ref{figure: sphere holes convergence} shows the mean value and the range of the relative error over the 50 experiments for different spatial step sizes $\Delta x$.
Second-order of convergence is observed for the mean value of the relative error.


\subsubsection{Reaction-diffusion systems}
This example considers a reaction-diffusion system, namely the Brusselator \cite{yang2004stable}, on a bumpy torus \cite{bumpytorus}. The system has the form
\begin{equation}\label{Brusselator}
  \begin{array}{l}
    u_t=\nu\Delta u+f(u,v),\\
    v_t=\mu\Delta v+g(u,v),
  \end{array}
\end{equation}
with
$$f(u,v) = a-(b+1)u+u^2v,\qquad g(u,v)=bu-u^2v,$$
where $a,b,\nu$ and $\mu$ are constants. The use of
the Semi-implicit 2-step Backward Differentiation Formula (SBDF2) scheme is recommended for the solution of reaction-diffusion equations when second-order centered finite difference
schemes are applied to the diffusion operator \cite{ruuth1995implicit}. Applying the proposed least-squares implicit RBF-CPM with the SBDF2 scheme yields
$$\left\{
  \begin{array}{l}
    \displaystyle 3U^{n+1}-4U^n+U^{n-1} = 2\Delta t(\nu W U^{n+1} + 2f(U^n,V^n)-f(U^{n-1},V^{n-1})),\\
    \displaystyle 3V^{n+1}-4V^n+V^{n-1} = 2\Delta t(\mu W V^{n+1} + 2f(U^n,V^n)-f(U^{n-1},V^{n-1})),\\
    \displaystyle U^{n+1} = PU^{n+1},\\
    \displaystyle V^{n+1} = PV^{n+1},\\
  \end{array}\right.
$$
where $W$ and $P$ are defined in equations~(\ref{Wu}) and (\ref{Pu}) respectively, and $U$ and $V$ are
the discretized solutions. A step of the implicit-explicit Euler discretization \cite{ascher1995implicit} is used to obtain starting
values for SBDF2.

\begin{figure}
    \centering
    \includegraphics[width=0.49\textwidth]{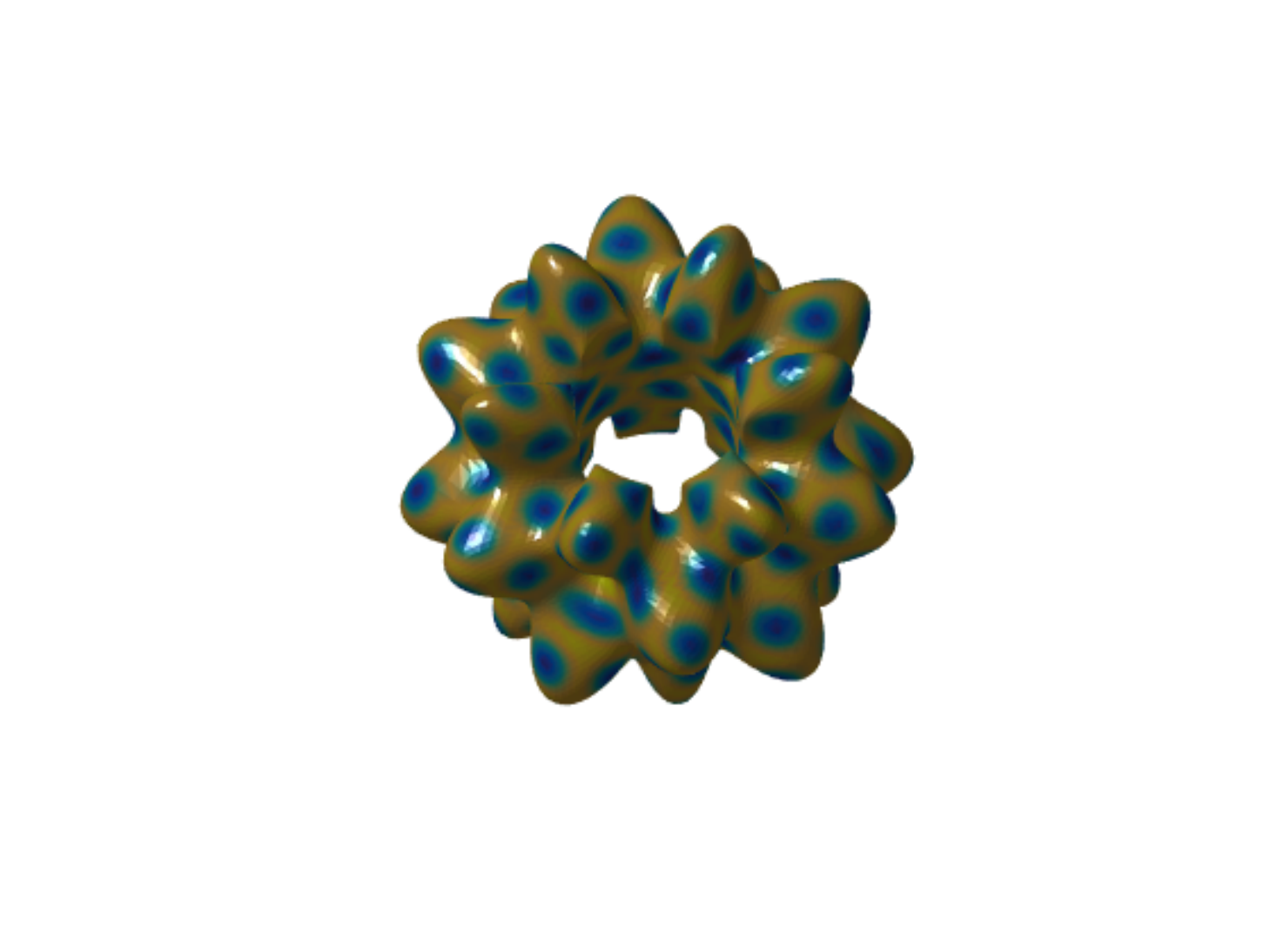}
    \includegraphics[width=0.49\textwidth]{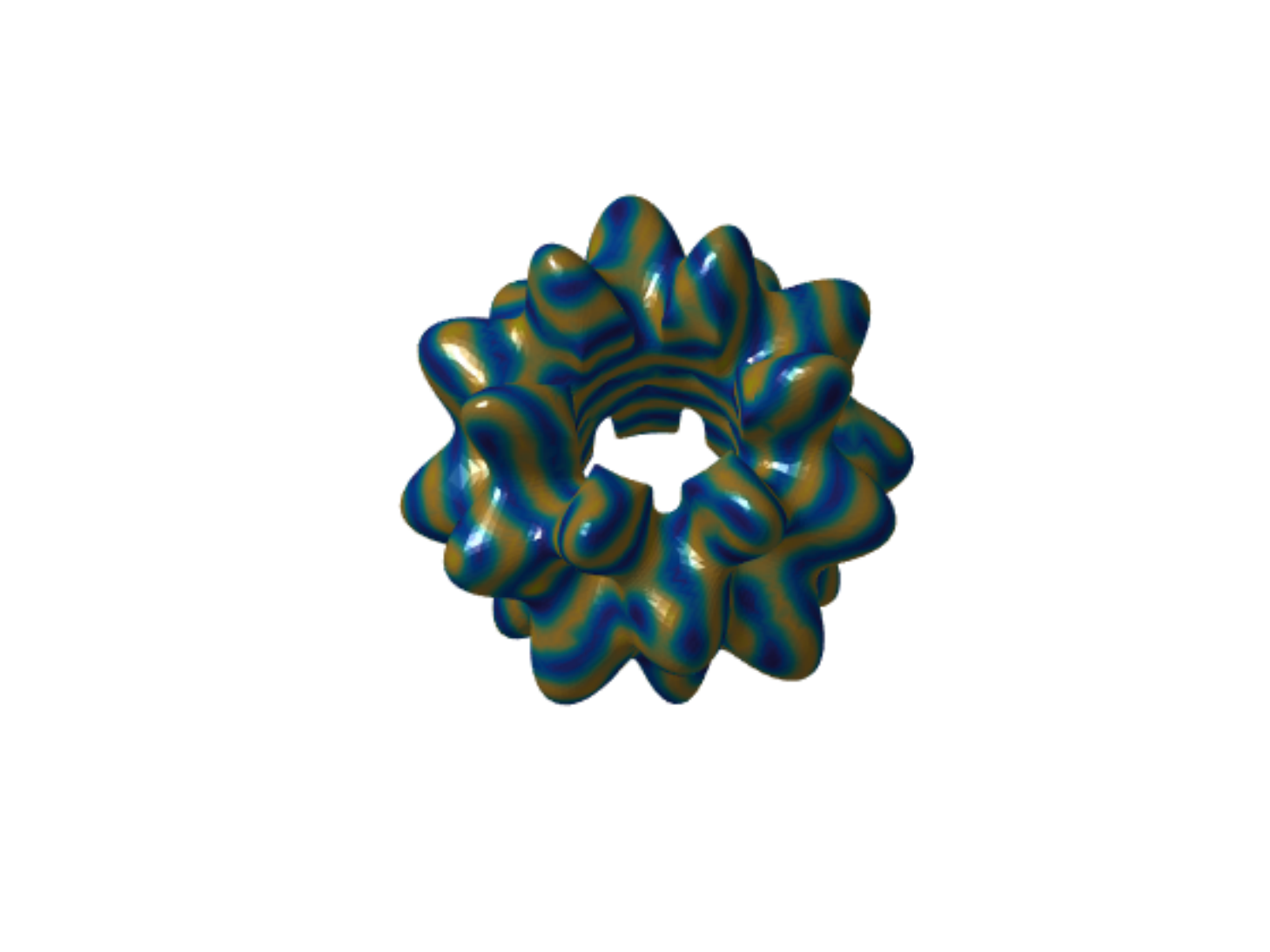}
    \caption{The solution of the reaction-diffusion system on a bumpy torus. The parameters chosen are $a=3$, $b=10.2$ and $\mu=10/900$.
    Finally, $\nu=5/900$ on the left figure and $\nu=3.8/900$ on the right figure.}
    \label{RDBumpyTorus}
\end{figure}

Two approximations of $U$ are shown in Figure~\ref{RDBumpyTorus}. In this numerical experiment, a random selection of
$1\%$ of points in the computational tube is removed.

\subsubsection{Cahn-Hilliard on a sphere} \label{sec:CHonSphere}
Our final example of PDEs on static surfaces considers the Cahn-Hilliard equation on the unit sphere. The equation
has the form \cite{gera2017cahn,gera2017stochastic}
\begin{equation}\label{equation: Cahn-Hilliard}
u_t - \frac{1}{P_e}\nabla_\Gamma\bigdot\left(\nu\nabla_\Gamma\left(\frac{\partial g}{\partial u}\right)\right) +
\frac{C_n^2}{P_e}\nabla_\Gamma\bigdot(\nu\nabla_\Gamma\Delta_\Gamma u)=0,
\end{equation}
where $C_n$ is the Cahn number, $P_e$ is the surface Peclet number, $\nu$ is the mobility and $g$ is a double well potential. For constant mobility $\nu=1$, equation~\eqref{equation: Cahn-Hilliard} can be split into two second-order equations \cite{gera2017cahn}
\begin{equation}\label{equation: Cahn-Hilliard system}
\begin{array}{c}
u_t-\frac{1}{P_e}\Delta_\Gamma \mu = 0,\\
\mu + C_n^2\Delta_\Gamma u = \frac{\partial g}{\partial u}(u),
\end{array}
\end{equation}
where $\mu$ is the chemical potential.

A number of discretizations are proposed in \cite{gera2017cahn}. We employ the first-order discretization scheme augmented with the constant-along-normal equations. Specifically,
\begin{equation}\label{CHDiscretization}
\left\{
\begin{array}{l}
U^{n+1}-U^n-\Delta t\frac{1}{P_e}WM^{n+1} = 0,\\
M^{n+1}+C_n^2WU^{n+1} = \frac{\partial g}{\partial u}(U^n),\\
U^{n+1} = PU^{n+1},\\
M^{n+1} = PM^{n+1},\\
\end{array}
\right.
\end{equation}
where $M$ and $U$ are the discretized approximations of $\mu$ and $u$, respectively, and $P$ and $W$ are defined as in \eqref{Pu} and \eqref{Wu}.

\begin{figure}
    \centering
    \includegraphics[width=0.49\textwidth]{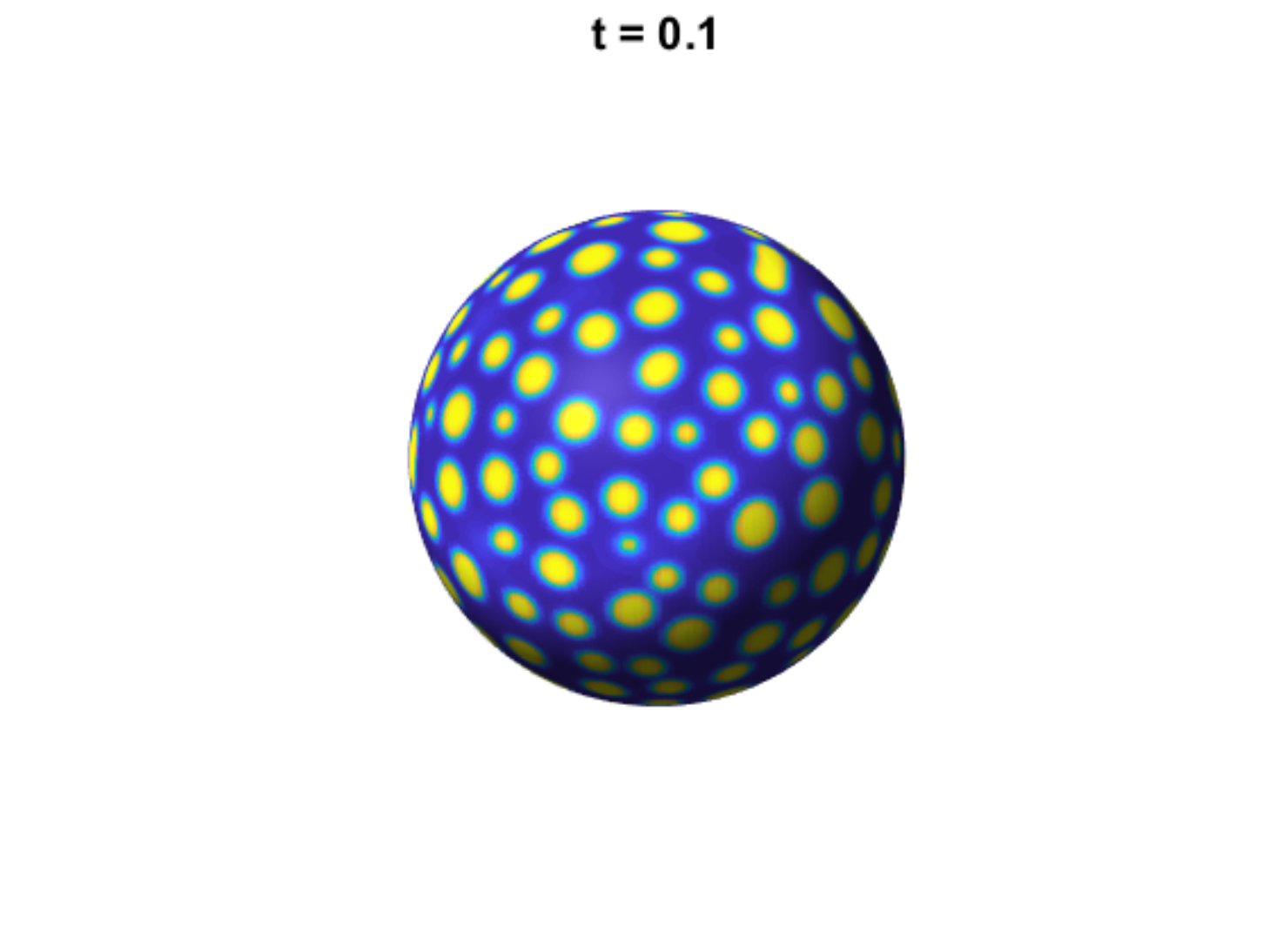}
    \includegraphics[width=0.49\textwidth]{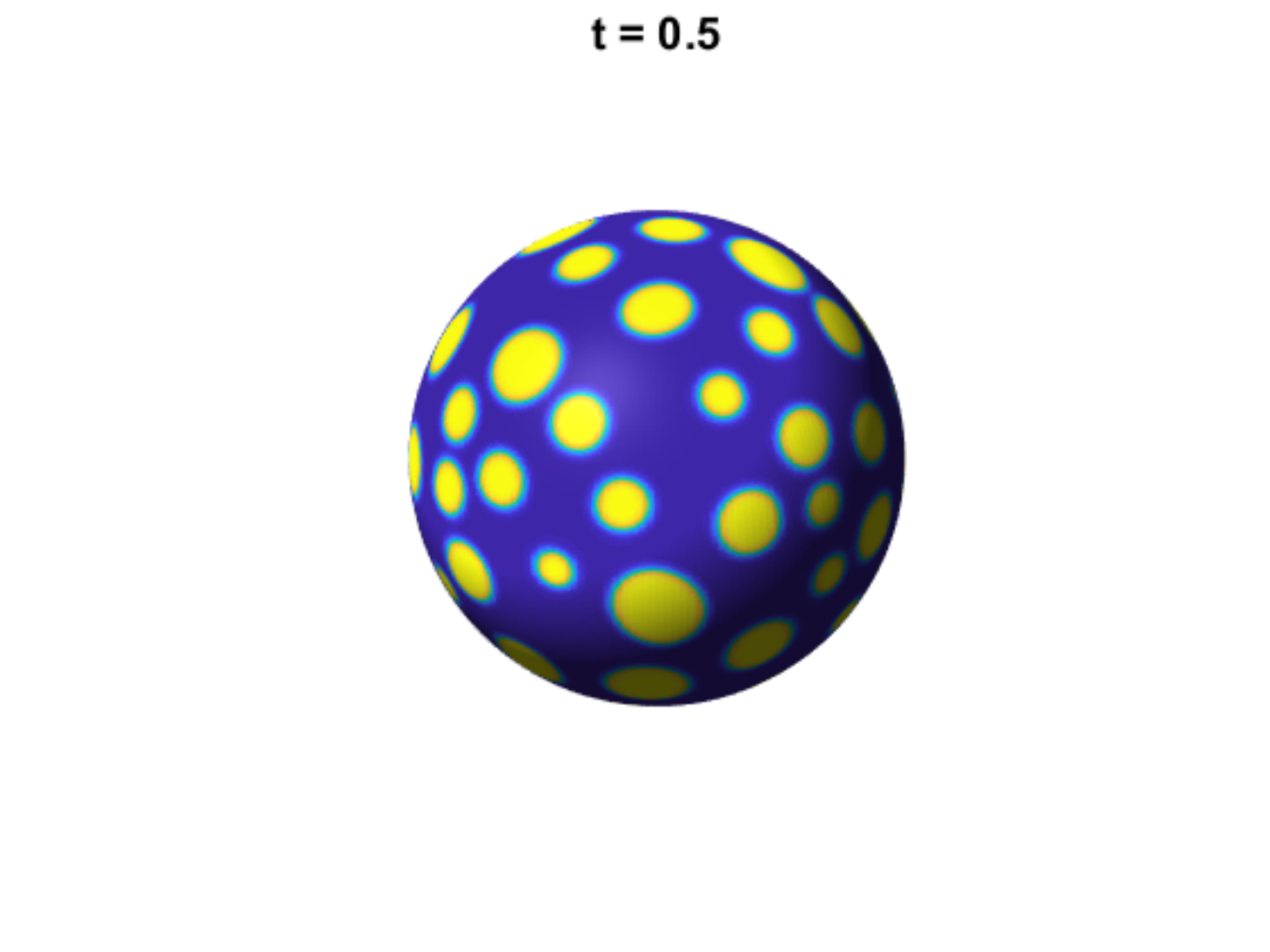}
    \includegraphics[width=0.49\textwidth]{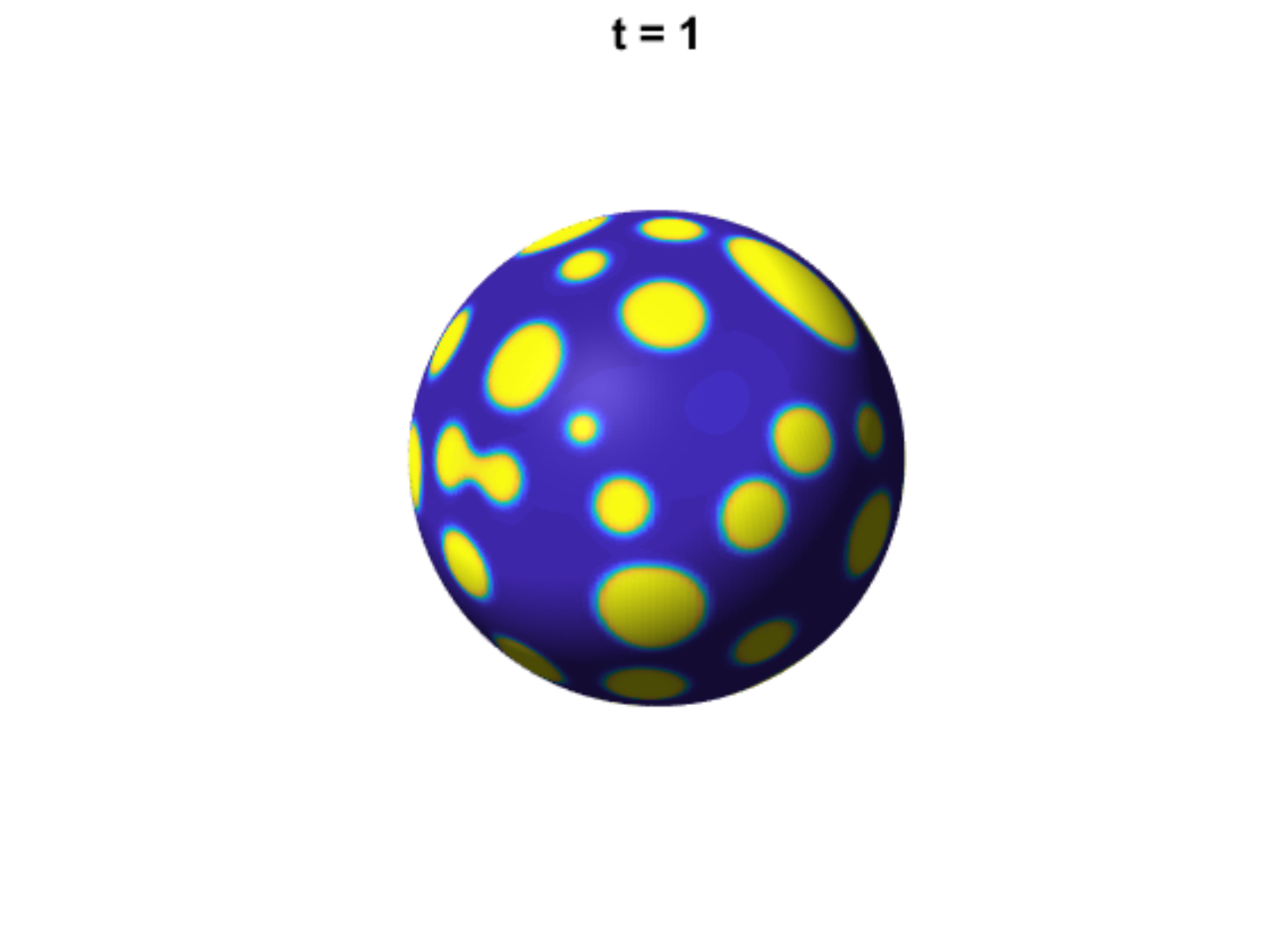}
    \includegraphics[width=0.49\textwidth]{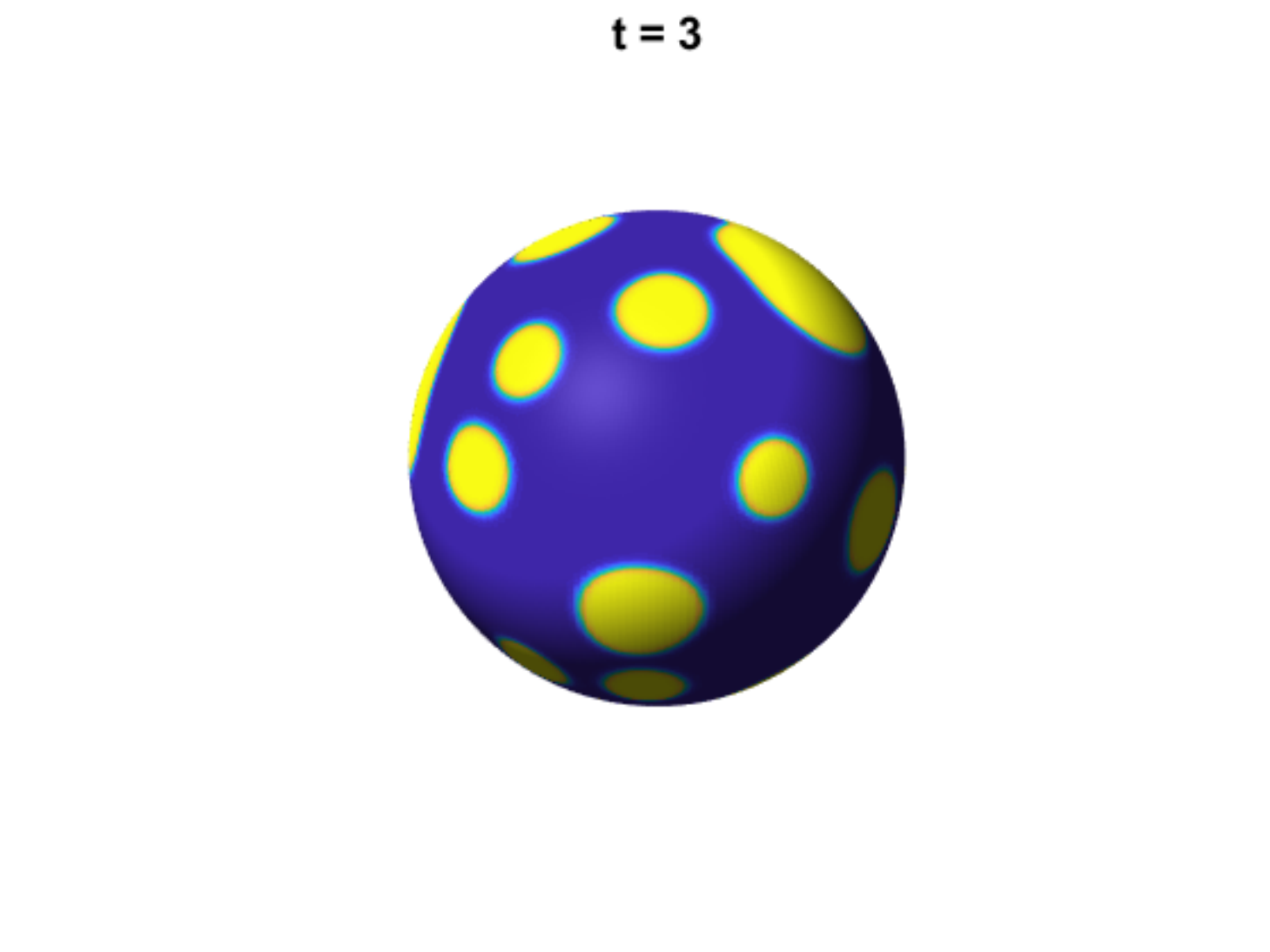}
    \caption{The solution of the Cahn-Hilliard equation on the unit sphere for various times $t$. Yellow and blue colors correspond to
    high and low concentrations respectively.}
    \label{solCahnHilSphere}
\end{figure}

Consider the unit sphere with an initial profile $u(\mathbf{x},0)$ defined as a random small
perturbation of magnitude 0.01 around the value 0.3. Further, select $\nu=1$, $P_e=1$, $C_n = 0.015$, and a double well potential function $g(u) = u^2(1-u)^2$.
For a spatial grid size $\Delta x = 0.0195$ and  a time step-size $\Delta t=10^{-4}$, chosen according to \cite{gera2017cahn}, we obtain the computed results
shown in Figure~\ref{solCahnHilSphere}. Due to the long computational time, a surface phase conservation is employed as described in \cite{gera2017stochastic}, by applying a correction factor $\beta$ at the solution at every time step, calculated as
$$\int_\Gamma \beta u^n\; dA = \int_\Gamma u^0\; dA,$$
where $u^0$ is the initial solution and $u^n$ the solution at time $n\Delta t$. The surface integrals are approximated 
using the singular values of the Jacobian of the closest point mapping \cite{kublik2016integration}. We use a centered difference scheme for the differentiation of the closest point mapping. The computational results demonstrate a similar behavior with the ones in \cite{gera2017stochastic}.

\section{A coupled method for solving PDEs on moving surfaces}\label{NumericalExamplesSection}
In this section, we couple the least-squares implicit RBF-CPM with the grid based particle method
to obtain a method for solving PDEs on {\it moving} surfaces.    Numerical experiments are provided to show the convergence of the
coupled method.

\subsection{A coupled method}
To begin, a Cartesian grid $Z = \{\mathbf{z}_j\}$ is constructed in the embedding space that contains the surface $\Gamma$. All the grid points in a neighborhood of $\Gamma$ of radius $\gamma$ are collected to form the computational tube $Z^0$. Using a closest point representation of the surface $cp_\Gamma^0$, the grid points in the computational tube are mapped to their closest points on the surface, $X^0 = cp_\Gamma^0(Z^0) := \{\mathbf{x}_j :\; \mathbf{x}_j = cp_\Gamma^0(\mathbf{z}_j), \text{ for }\mathbf{z}_j\in Z^0\}$. Given an initial solution $u^0$ on the surface $\Gamma$, a constant-along-normal extension is defined on $Z^0$ using the closest point mapping $cp_\Gamma^0$; $U^0_Z = u^0(cp_\Gamma^0(\mathbf{z}_j)) = u^0(\mathbf{x}_j)$. We time step from $t = n\Delta t$ to $t = (n+1)\Delta t$, $n = 0,1,2,...$ according to the following combination of the GBPM and the RBF-CPM:

\begin{enumerate}
\item \textbf{Motion:} Evolve the points $X^n$, according to the desired motion law to yield $X^*$.
\item For all $\mathbf{z}_j\in Z^n$ and $\mathbf{x}_j\in X^n = cp_\Gamma^n(Z^n)$,
	\begin{enumerate}
	\item \textbf{Resampling:} Find the new closest point $\mathbf{x}^{n+1}_j$ of the grid point $\mathbf{z}_j$, using a local least-squares polynomial reconstruction on the moved surface points $X^*$.
	\item \textbf{RBF-FD calculation:} Calculate the RBF-FD weights $\mathbf{w}_j$ in \eqref{eq:RBF-FD weights} for $\mathbf{x}_j$, using the grid points $Z^n$. Also, if $\mathbf{x}^{n+1}_j$ can be found and its distance from $\mathbf{z}_j$ is smaller than the tube radius $\gamma$, calculate the RBF-FD weights $\mathbf{w}_j^{n+1}$ in \eqref{eq:RBF-FD weights} for $\mathbf{x}^{n+1}_j$, using the grid points $Z^n$.
	\end{enumerate}
\item For all neighboring grid points $\mathbf{z}_k$ of the computational tube $Z^n$ repeat steps (a) and (b) of step 2, using $\mathbf{z}_k$ in place of $\mathbf{z}_j$ everywhere.
\item \textbf{Point deactivation:} Deactivate the grid points $Z^{n+1} = \{\mathbf{z}_j :\; \mathbf{x}^{n+1}_j = cp_\Gamma^{n+1}(\mathbf{z}_j)\}$ and their corresponding closest points in $X^{n+1} = \{\mathbf{x}^{n+1}_j\}$, if their distance is larger than $\gamma$.
\item \textbf{PDE solution:} Solve the embedding PDE and get the solution $U^*_Z$ at the grid points $Z^n$.
\item \textbf{Interpolation:} Interpolate the solution $U_Z^*$ from $Z^n$ to the surface points $\{\mathbf{x}^{n+1}_j\} := X^{n+1} = cp_\Gamma^{n+1}(Z^{n+1})$ to get the solution $U^{n+1}_Z$. 
\end{enumerate}

Note that in the algorithm above, all the grid points $\mathbf{z}_j$ that can successfully be mapped to their closest points $\mathbf{x}^{n+1}_j$ and are within the tube radius $\gamma$ in steps 2 and 3 are contained to the computational tube $Z^{n+1}$. This mapping is denoted as $cp_\Gamma^{n+1}$ in the algorithm above. Also, step 3 (update of the computational tube) might not be necessary for each time step $\Delta t$, thus an extra condition can be applied to identify the candidates $\mathbf{z}_k$ neighboring $Z^n$ for which the step is required:
$$\min_j\|\mathbf{z}_k-\mathbf{x}^*_j\| \leq c_1\gamma,$$
where $\mathbf{x}^*_j\in X^*$ are the points on the moved surface in step 1 and $c_1\geq 1$ is a constant. If this condition is satisfied, then the application of step 3 for $\mathbf{z}_k$ is necessary.
For more details on the RBF-FD calculation, we refer to Section~\ref{ImplicitCPMDescription} and \cite{petras2018explicit}. Additional information on the resampling step of the GBPM can be found in \cite{leung2009grid}.

We apply a similar concept to \cite{petras2016pdes} for the calculation of the computational tube radius $\gamma$.
For a static surface and an $m$-point RBF-FD stencil, a computational tube radius $\tilde{\gamma}$ may be determined using the Gauss circle problem \cite{petras2018explicit}.  Then, for a moving surface, a safe choice for the tube
radius is
\begin{equation} \label{eq:safegamma2}
\gamma = \tilde{\gamma}+\Delta t \bigdot v_n^{max},
\end{equation}
where $v_n^{max}$ is an upper bound on the speed of a footpoint in the normal direction at time step $n$.
In practice, this is too conservative and we select $\gamma$ adaptively by setting $\gamma=\tilde{\gamma}$ and checking for violations of the tube radius. These violations of the tube radius are infrequent, and when they arise we recompute using
a larger computational tube (e.g., using (\ref{eq:safegamma2})), while keeping the same number of points $m$ in the RBF-FD stencil.

\subsection{Numerical experiments}\label{sec:Moving numerical examples}
In this section, we present several numerical experiments for the solution of PDEs on moving surfaces.
Unless stated otherwise, we approximate the solution of the advection-diffusion PDE (\ref{AdvDifPDE2}) with a diffusivity parameter $\mathcal{D}=1$. Due to
our use of a closest point representation, the normal derivative vanishes in equation~(\ref{AdvDifPDE2}): $\partial u/\partial n=0$.
Implicit-explicit Euler time stepping is chosen \cite{ascher1995implicit}. This yields
\begin{equation}\label{AdvDifDiscretization}
\left\{\begin{array}{l}
\displaystyle U^{n+1} -\Delta tW U^{n+1} = U^n+\Delta t(f(U^n) -V^n\kappa^n U^n+\sum_k^d D_k(U^nT^n_k)),\\
\displaystyle U^{n+1} = PU^{n+1},\\
\end{array}
\right.
\end{equation}
where $U^n$ is the discretized solution at time $n\Delta t$, $P$ and $W$ are defined as in \eqref{Pu} and \eqref{Wu}, $T^n_k$ is the $k$-th component of the tangential velocity $\mathbf{T}$ of the surface, $d$ is the dimension of the embedding space and $D_k$ is the $k$-th first derivative approximation using RBF-FD, i.e. $k=1$ corresponds to the first derivative in $x$, for which the RBF-FD coefficients for each $\mathbf{x}_j$ are calculated using
$$\frac{\partial \tilde{u}}{\partial x}(\mathbf{x}_j)\approx\mathbf{w}_j\tilde{u}(Z_j),$$
where $Z_j$ is the set of the $m$ closest grid nodes to $\mathbf{x}_j$. The matrix $D_1$ is constructed using $\mathbf{w}_j$ from the equation above, as described in Section~\ref{ImplicitCPMDescription}. The PHS RBF of order 7 is used with augmented polynomial basis that spans $\mathbb{P}_3$ for the calculation of all the matrices $P$, $D_k$ and $W$. 
The system~(\ref{AdvDifDiscretization}) is solved in a least-squares sense in the same manner
as described previously in Section~\ref{ImplicitCPMDescription}. 

\subsubsection{Diffusion on an expanding circle}
In our first example, we investigate the convergence of our method in two dimensions.
Following \cite{petras2016pdes}, we consider the homogeneous PDE (\ref{AdvDifPDE2})
on a circle centered at the origin with an initial radius $r_0=0.75$. A constant normal velocity $\mathbf{v}=5\mathbf{n}$ is imposed.
With these choices, the exact solution is
$$u(\theta,t)=e^{4/(5r(t))}\frac{\cos\theta\;\sin\theta}{r(t)},$$
where $r(t)=r_0+5t$ is the radius of the circle at time $t$.

\begin{figure}
    \centering
    \includegraphics[width=0.6\textwidth]{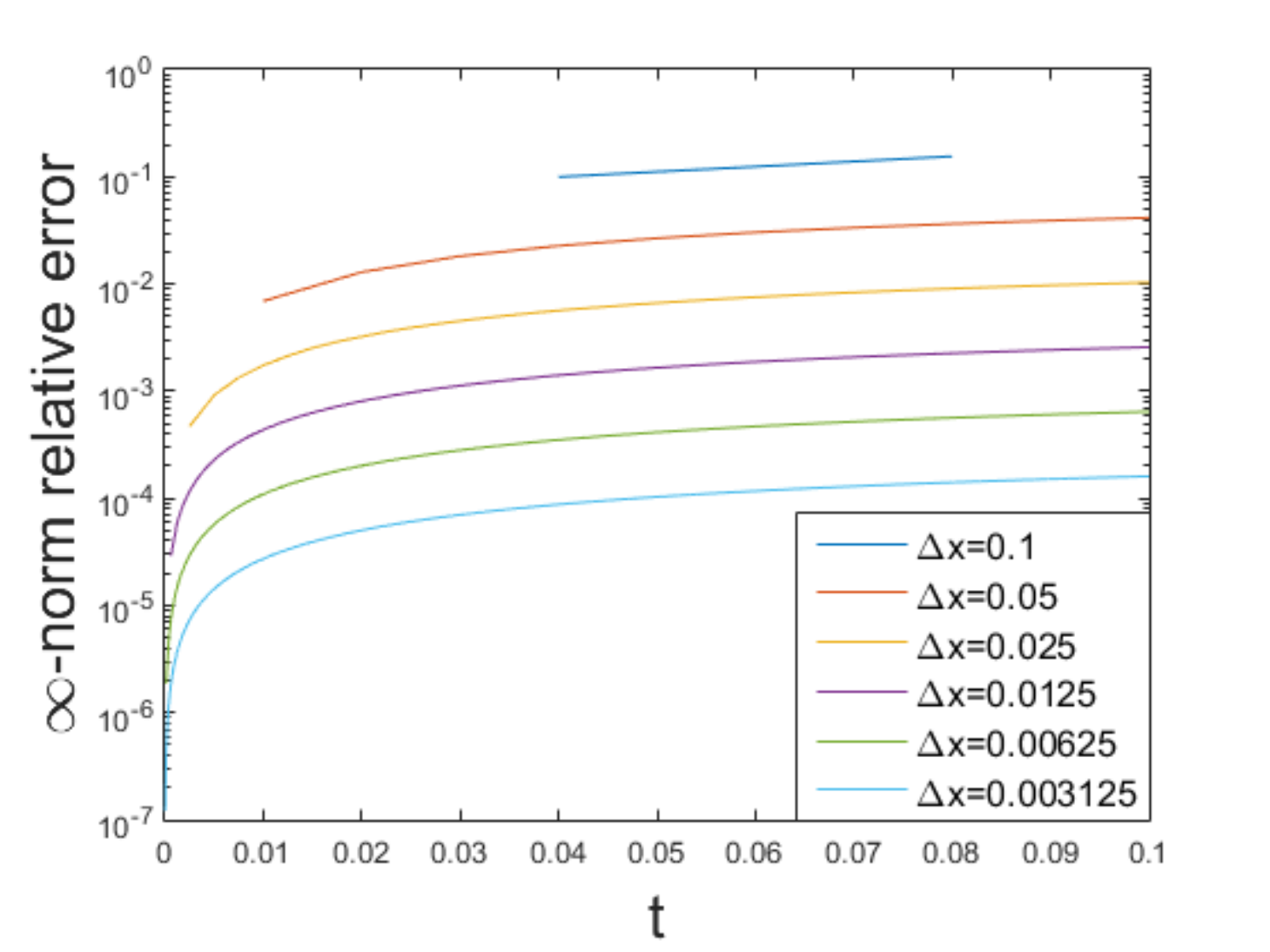}
    \caption{The relative error of the numerical solution of the diffusion model
     on the expanding circle over time. Errors are computed using the analytical solution.}
    \label{errorMovingCircle}
\end{figure}

\begin{table}
\centering
  \begin{tabular}{|l|cc|cc|}
    \hline
    $\Delta x$ & $t =$ 0.04 & e.o.c. & $t =$ 0.08 & e.o.c.  \\ \hline
    0.1 & 9.96$\times 10^{-2}$ & - & 1.55$\times 10^{-1}$ & - \\
    0.05 & 2.27$\times 10^{-2}$ & 2.13 & 3.63$\times 10^{-2}$ & 2.09\\
    0.025 & 5.60$\times 10^{-3}$ & 2.02 & 9.00$\times 10^{-3}$ & 2.01\\
    0.0125 & 1.40$\times 10^{-3}$ & 2.01 & 2.30$\times 10^{-3}$ & 2.00 \\
    0.00625 & 3.51$\times 10^{-4}$ & 2.00 & 5.63$\times 10^{-4}$ & 2.00\\
    0.003125 & 8.75$\times 10^{-5}$ & 2.00 & 1.41$\times 10^{-4}$ & 2.00\\
    \hline
  \end{tabular}
  \caption{Relative errors as measured in the infinity norm and the estimated order of convergence (e.o.c.) at various times $t$.}
  \label{errorMovingCircleTable}
\end{table}

Using our discretization~(\ref{AdvDifDiscretization}),  the relative
error of the approximate solution (i.e., compared to the known, exact solution) is computed over time  for different grid sizes $\Delta x$; see Figure~\ref{errorMovingCircle}.
Table~\ref{errorMovingCircleTable} shows the relative errors and the corresponding convergence rates
for selected times and grid sizes. Second-order convergence is observed
for a time step-size of $\Delta t=4\Delta x^2$.

\subsubsection{Advection-diffusion on an oscillating sphere}
Next, we consider advection-diffusion on an oscillating ellipsoid.
Following the example in \cite{petras2016pdes}, we evolve the  non-homogeneous PDE (\ref{AdvDifPDE2}) starting from the unit sphere. The velocity of the ellipsoid is a function of the $x$-coordinate:
$$\mathbf{v} = \frac{a'(t)}{2a(t)}(x_1,0,0),$$
where
$$a(t) = 1+\sin(2t).$$
For these choices, the exact shape of the surface is simply
$$\mathbf{x}(\theta,\phi,t) = (\sqrt{a(t)}\cos(\theta)\cos(\phi),\sin(\theta)\cos(\phi),\sin(\phi)),$$
where $\theta$ is the azimuth angle and $\phi$ is the elevation angle.
The forcing term of the PDE is chosen to be
$$f = u(\mathbf{x},t) \bigdot \Big(-6+\frac{a'(t)}{a(t)}\Big(1-\frac{x_1^2}{2N}\Big)+\frac{1+5a(t)+2a^2(t)}{N}-
\frac{1+a(t)}{N^2}(x_1^2+a^3(t)(x_2^2+x_3^2))\Big),$$
with
$$N = x_1^2+a^2(t)(x_2^2+x_3^2).$$
For this choice of $f$, the solution of the PDE is
$$u(\mathbf{x},t) = e^{-6t}x_1x_2,$$
for all times $t\ge0$.

\begin{figure}
    \centering
    \includegraphics[width=0.6\textwidth]{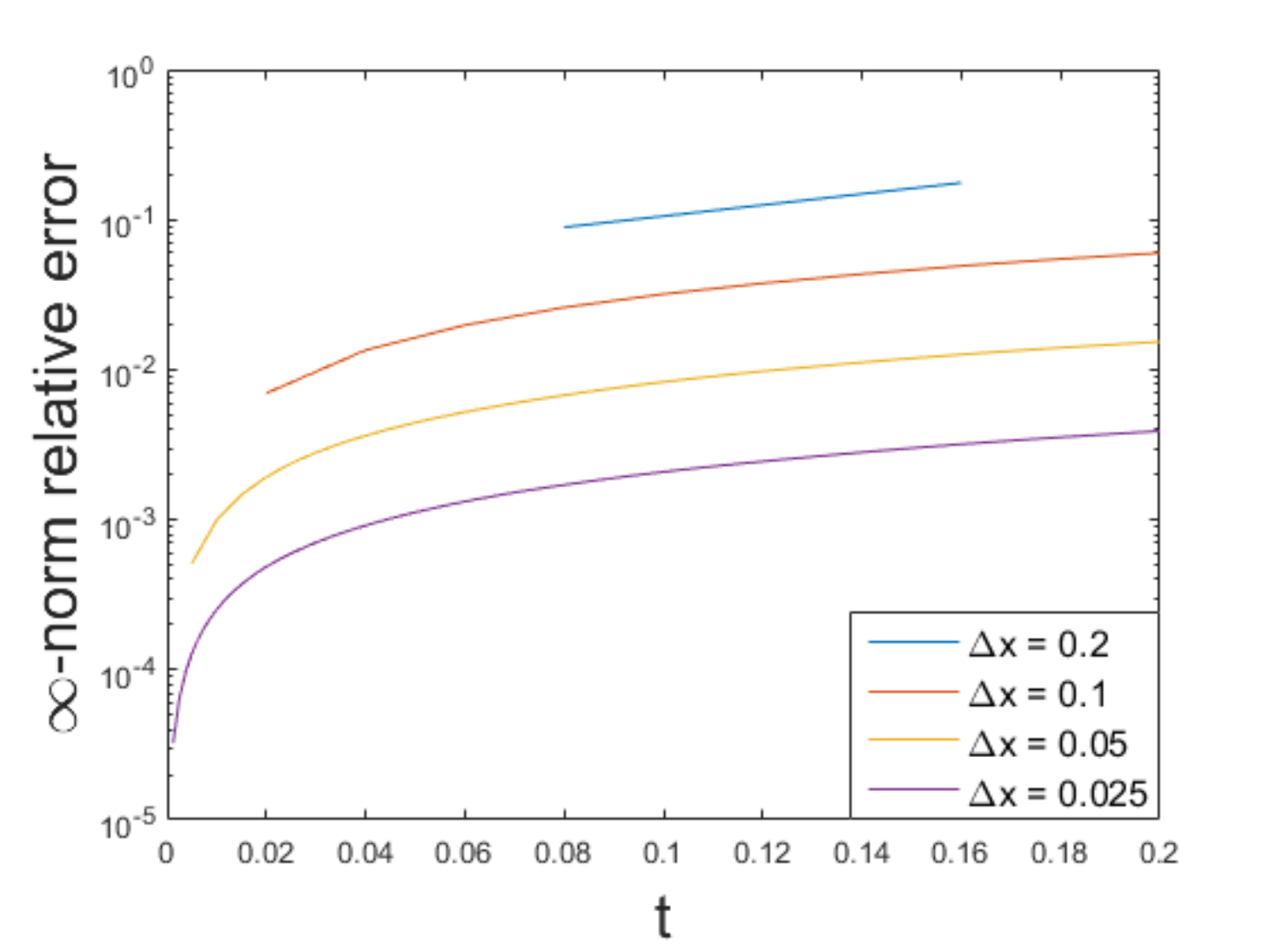}
    \caption{The relative error of the numerical solution of the diffusion model
     on the oscillating sphere over time. Errors are computed using the analytical solution.}
    \label{errorMovingSphere}
\end{figure}

\begin{table}
\centering
  \begin{tabular}{|l|cc|cc|}
    \hline
    $\Delta x$ & $t =$ 0.08 & e.o.c. & $t =$ 0.16 & e.o.c.  \\ \hline
    0.2 & 8.89$\times 10^{-2}$ & - & 1.76$\times 10^{-1}$ & - \\
    0.1 & 2.6$\times 10^{-2}$ & 1.77 & 4.91$\times 10^{-2}$ & 1.84\\
    0.05 & 6.80$\times 10^{-3}$ & 1.94 & 1.26$\times 10^{-2}$ & 1.96\\
    0.025 & 1.70$\times 10^{-3}$ & 1.99 & 3.20$\times 10^{-3}$ & 1.99 \\
    \hline
  \end{tabular}
  \caption{Relative errors as measured in the infinity norm and the estimated order of convergence (e.o.c.) at various times $t$.}
  \label{errorMovingSphereTable}
\end{table}

Using the implicit-explicit Euler method~(\ref{AdvDifDiscretization}) with $\Delta t = 2\Delta x^2$,
we compute and plot (in Figure~\ref{errorMovingSphere})  the
$\infty$-norm relative  error of the approximation  over time for different spatial grid sizes $\Delta x$.
Errors and the corresponding convergence rates at selected times $t$ are presented in Table~\ref{errorMovingSphereTable}.
Second-order convergence in $\Delta x$ is observed.

\subsubsection{A cross-diffusion reaction-diffusion system}
This example considers a reaction-diffusion system with cross-diffusion terms \cite{madzvamuse2014exhibiting}.
The system of equations has the form
\begin{equation} \label{eq:rdmodel}
    \begin{array}{l}
        \displaystyle\frac{Du}{Dt} = \Delta_{\Gamma}u+d_w\Delta_{\Gamma}w-u\nabla_{\Gamma}\bigdot\mathbf{v}+f_1(u,w),\\
        \\
        \displaystyle\frac{Dw}{Dt} = \mathcal{D}\Delta_{\Gamma}w+d_w\Delta_{\Gamma}u-w\nabla_{\Gamma}\bigdot\mathbf{v}+f_2(u,w),\\
  \end{array}
\end{equation}
where $\mathcal{D},d_u$ and $d_w>0$ are scalar quantities. The coupling functions that we consider are
$$f_1(u,w) = 200(0.1-u+u^2w),\qquad f_2(u,w) = 200(0.9-u^2w).$$
We select a velocity that is purely in the normal direction,
$$\mathbf{v}=(0.01\kappa+0.4u)\mathbf{n},$$
where $\kappa$ is the mean curvature and $\mathbf{n}$ is the unit normal vector. In our experiment, we
start from a torus, and set the initial $u$ and $v$ to be small random
perturbations around 0.5.    Our remaining parameters are chosen to be
 $\mathcal{D} = 10$ and $d_u = d_w=1$.
An implicit-explicit Euler time discretization is applied, specifically,
$$\left\{\begin{array}{l}
        U^{n+1} -\Delta tW U^{n+1}= U^n+\Delta t(f_1(U^n,W^n)
-V^n\kappa^n U^n-\mathbf{D}\bigdot(U^n\mathbf{T}^n)+d_wWW^n),\\
        W^{n+1} -\Delta t\mathcal{D}W W^{n+1} = W^n+\Delta t(f_2(U^n,W^n)
-V^n\kappa^n W^n-\mathbf{D}\bigdot(W^n\mathbf{T}^n)+d_uWU^n),\\
\displaystyle U^{n+1} = PU^{n+1}, \\
\displaystyle W^{n+1} = PW^{n+1}, \\
  \end{array}\right.$$
with a time step-size $\Delta t=0.2\Delta x^2$ and a spatial grid size $\Delta x=0.05$, where $W$ (not to be confused with the computational solution $W^n$), $P$ and $\mathbf{D}$ are defined as at the beginning of Section~\ref{sec:Moving numerical examples}.

\begin{figure}
    \centering
    \includegraphics[width=0.49\textwidth]{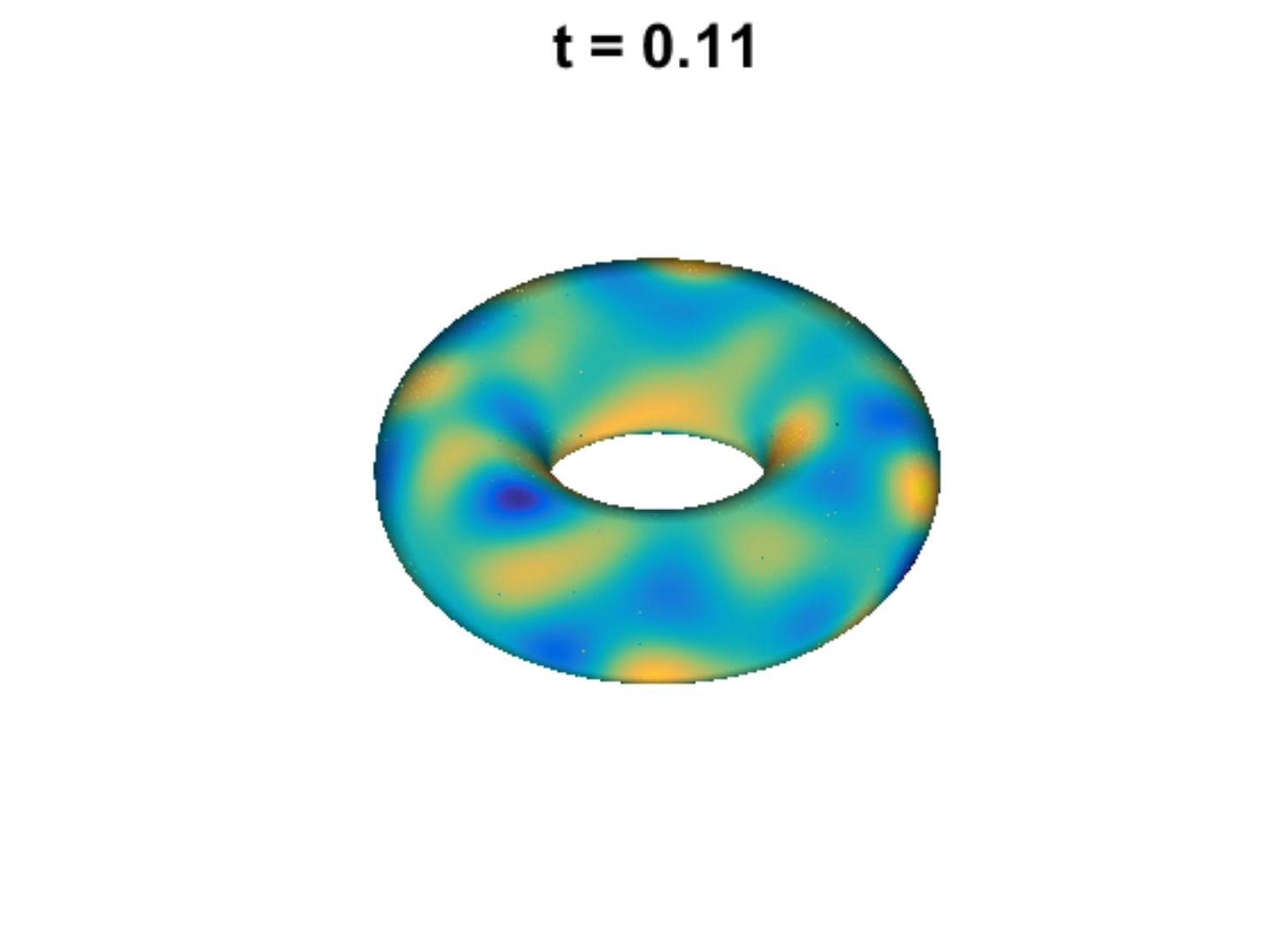}
    \includegraphics[width=0.49\textwidth]{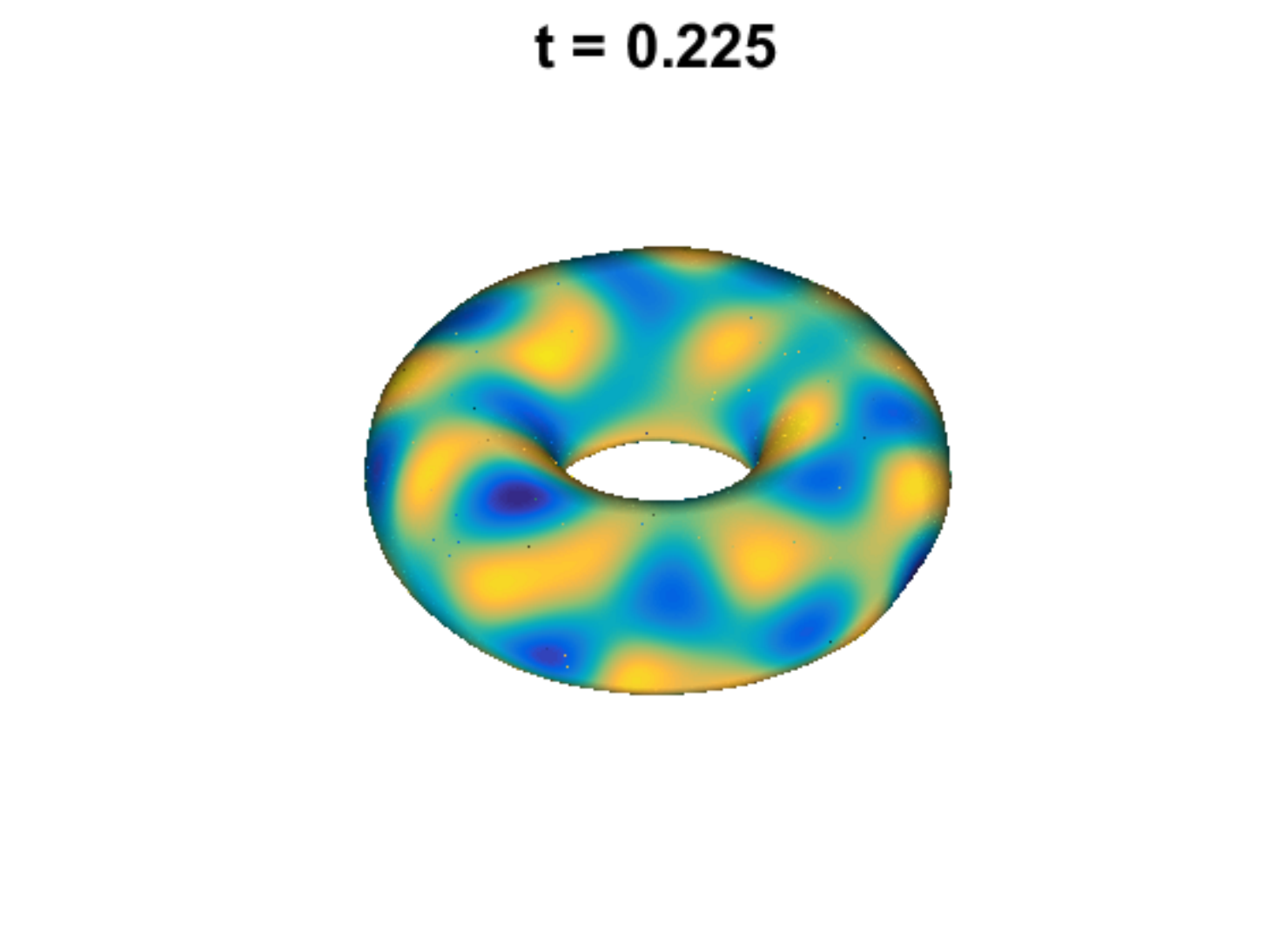}
    \includegraphics[width=0.49\textwidth]{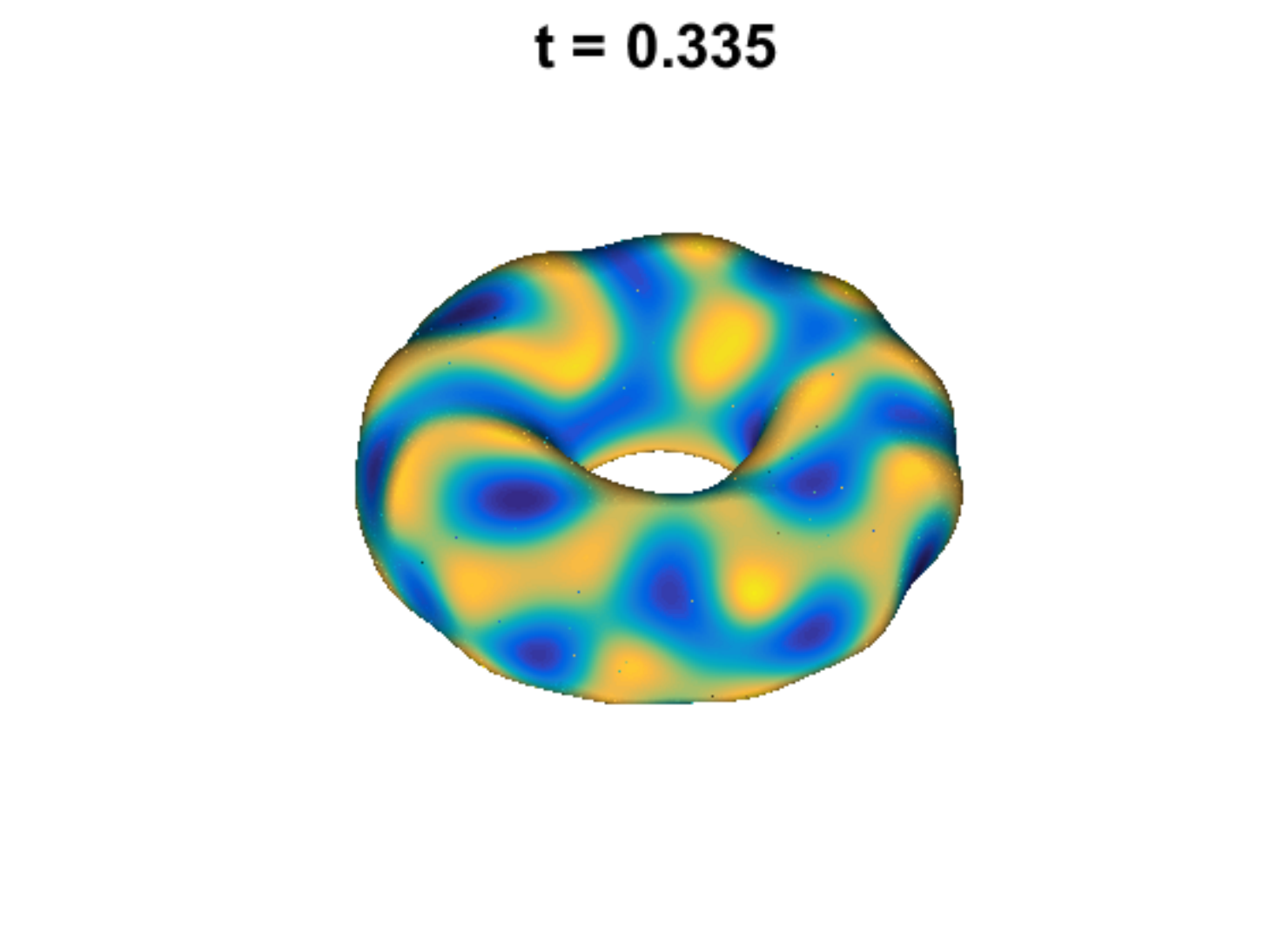}
    \includegraphics[width=0.49\textwidth]{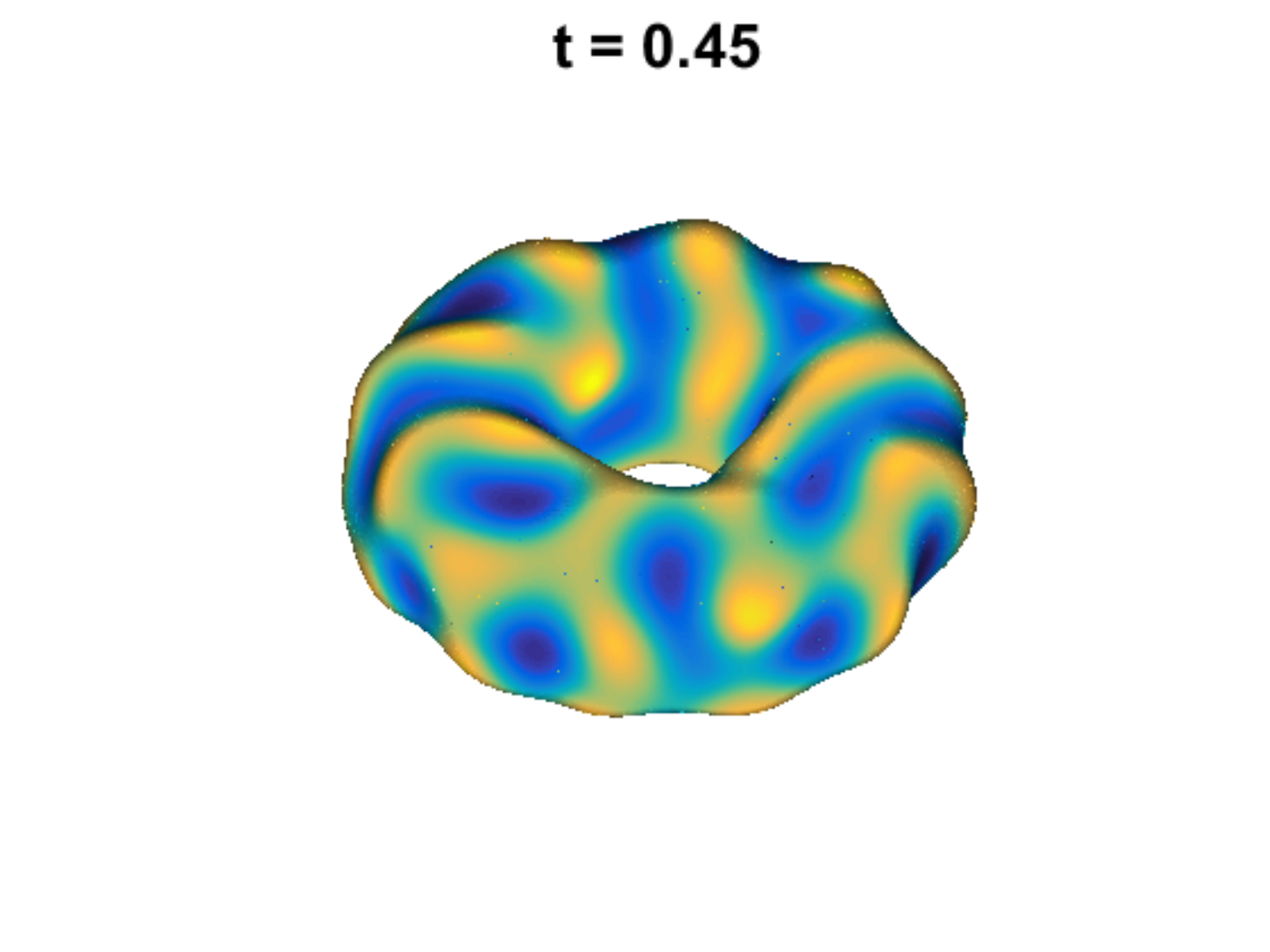}
    \caption{A visualization of the numerical solution $u$ of the cross-diffusion
reaction-diffusion system~(\ref{eq:rdmodel}) at selected times $t$.}
    \label{solRDTorus}
\end{figure}

Figure~\ref{solRDTorus} displays the results of our coupled method for selected times.
We observe that the surface remains smooth, and extends outwards in regions where the solution $U$ is large (areas shown in yellow).
See \cite{madzvamuse2014exhibiting} for further results for this model on evolving surfaces.

\subsubsection{Cahn-Hilliard equation on an ellipsoid}
Our final example considers the homogeneous Cahn-Hilliard PDE~(\ref{CahnHilPDE}). An initial ellipsoid centered at the origin with
semi-major axes $r_x=1.2$, $r_y=0.7$, and $r_z=0.7$ is evolved in the normal direction according to the velocity
$$\mathbf{v}=(0.01\kappa+0.4u)\mathbf{n},$$
where $\kappa$ is the mean curvature and $\mathbf{n}$ is the unit outward normal vector. We explore the effect of the surface motion to the Cahn-Hilliard solution, by using the same parameter setup as in Section~\ref{sec:CHonSphere}.

\begin{figure}
    \centering
    \includegraphics[width=0.49\textwidth]{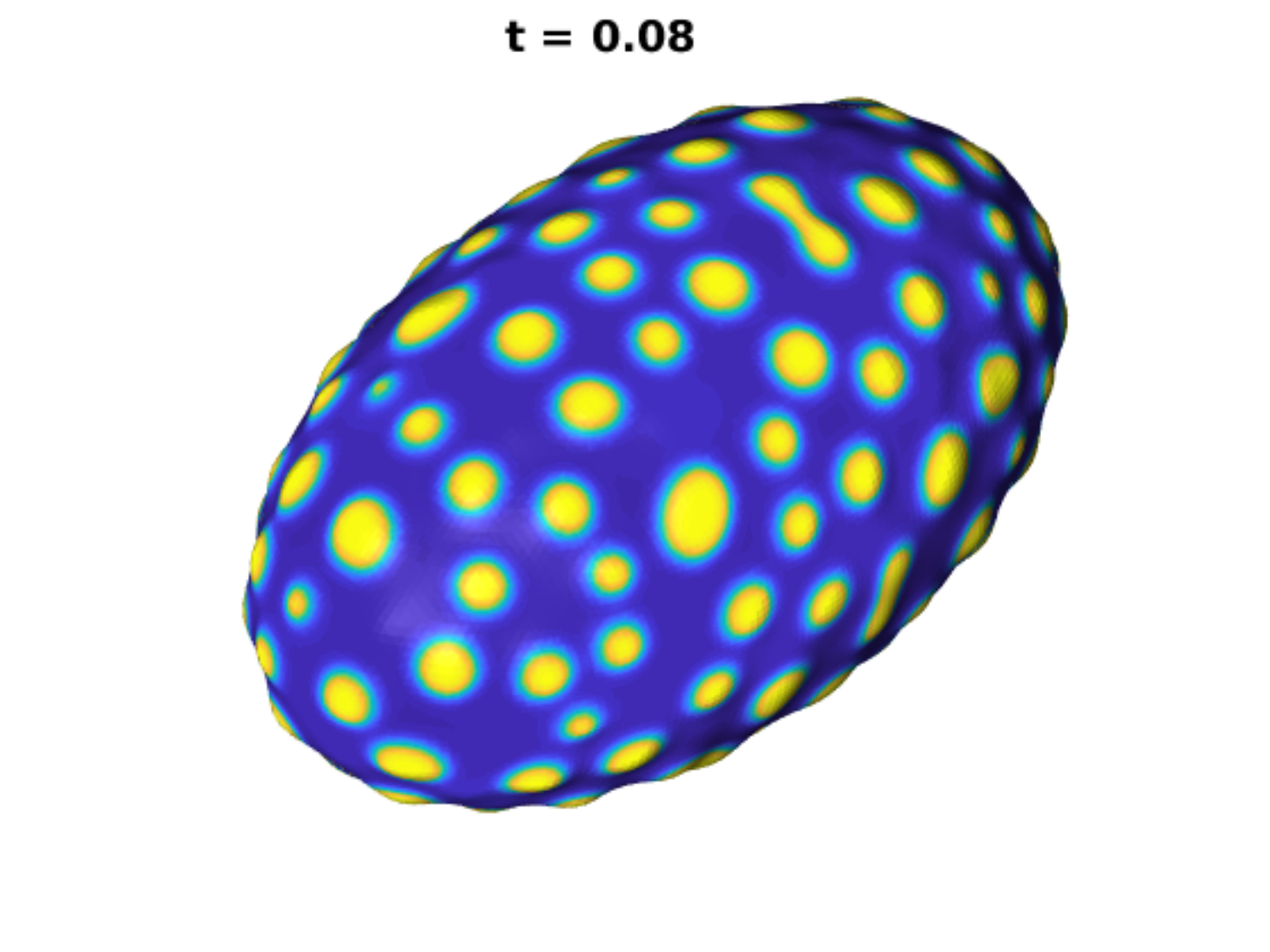}
    \includegraphics[width=0.49\textwidth]{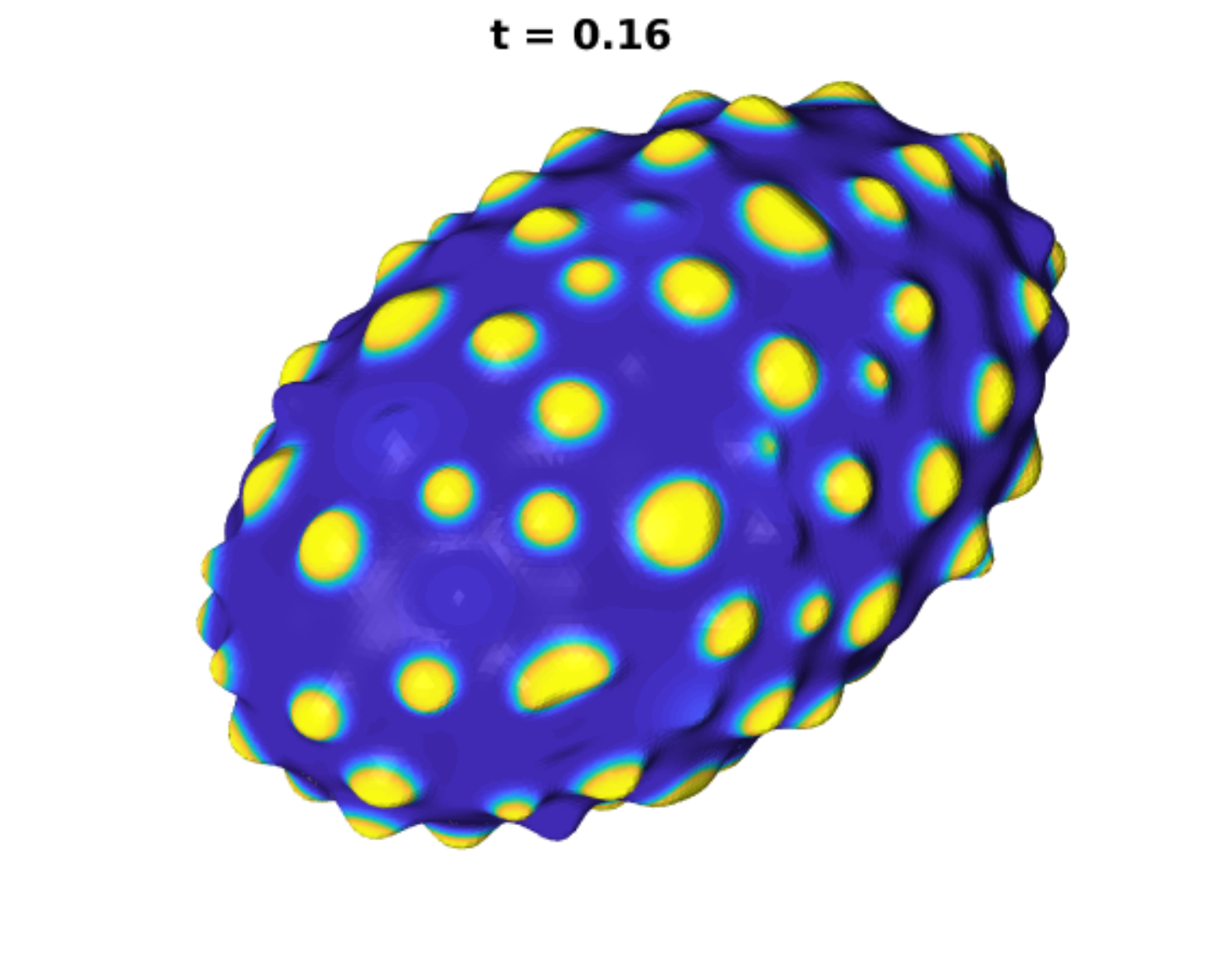}
    \caption{A visualization of the numerical solution $U^n$ of the Cahn-Hilliard PDE~(\ref{CahnHilPDE})
at selected times $t=n\Delta t$.}
    \label{solCHEllipsoid}
\end{figure}

The discretization is carried out according to the procedure
described in Section~\ref{sec:CHonSphere} with a spatial step-size $\Delta x = 0.02$ and a time step-size $\Delta t = 10^{-4}$. Specifically,
$$\left\{
\begin{array}{l}
U^{n+1}-U^n-\Delta t\left(\frac{1}{P_e}WM^{n+1} - V^n\kappa^nU^n\right) = 0,\\
M^{n+1}+C_n^2WU^{n+1} = \frac{\partial g}{\partial u}(U^n),\\
U^{n+1} = PU^{n+1},\\
M^{n+1} = PM^{n+1},\\
\end{array}
\right.$$
since we only consider a surface evolution in the normal direction. 
Selecting the initial profile $U^0$ to consist of small random perturbations
around 0.3, we obtain the numerical results displayed in Figure~\ref{solCHEllipsoid}.
We observe that the surface grows outwards at the high concentration regions of the solution $U$, shown in yellow.

\section{Summary}\label{SummarySection}
In this paper, we develop a novel least-squares approach to stabilize the closest point method using RBF-FD. In our method, the constant
normal extension of the approximate solution in the embedding
space is imposed using a separate equation.  The PDE and the constant extension in the
normal direction system is solved as a least-squares formulation. The least-squares
method provides the flexibility to use a variety of computational tubes around the surface,
including ones that have holes. Numerical examples illustrate the convergence of the method
on computational tubes missing 1\% and 5\% of the grid points (and their corresponding closest points
on the surface) that lie within.

Furthermore, a simple coupling of the least-squares implicit RBF-CPM and the grid based particle method is proposed to approximate the solution of the
conservation law described in Section~\ref{PDEsSection}. Numerical results show second-order convergence
for time step-sizes $\Delta t = \mathcal{O}(\Delta x^2)$. In addition, the coupled method is tested on strongly
coupled systems including a cross-diffusion reaction-diffusion model and the Cahn-Hilliard
equation.

Much work needs to be done to test the coupled method, including examples of PDE
models on moving open surfaces and on surfaces with topological changes. An interesting
application would employ the conservation law on two circles expanding in the normal
direction and merging into one curve \cite{petras2016pdes}.
The least-squares implicit RBF-CPM should be capable of dealing with the discontinuity that arises from the merging of
the two circles.

Moreover, the least-squares implicit closest point method using RBF-FD should
be tested for the case of adaptive computational tubes surrounding the surface. Adaptivity
is often required to capture fine features of a surface, i.e. areas with high curvature. The
grid based particle method contains an optional adaptivity step for such occurrences \cite{leung2009grid}.
It is expected that the RBF-FD implementation of the closest point method
will provide a simple way of developing schemes for adaptive tubes.

Another interesting direction for solving PDEs on moving surfaces employs meshfree
techniques that use the closest point concept. Such methods have been developed for
approximating the solution of surface PDEs \cite{piret2012orthogonal,cheung2015localized}, however no one until this paper has used
RBF approximations to solve PDEs on moving surfaces.

\section*{Acknowledgements}
The first and fourth authors were partially supported by an NSERC Canada grant (RGPIN 2016-04361). The first author was partially supported by the Basque Government through the BERC 2018-2021 program and by Spanish Ministry of Science, Innovation and Universities through The Agencia Estatal de Investigacion (AEI) BCAM Severo Ochoa excellence accreditation SEV-2017-0718 and through project MTM2015-69992-R BELEMET. This work was partially supported by a Hong Kong Research Grant Council GRF Grant, and a Hong Kong Baptist University FRG Grant. This research was enabled in part by support provided by WestGrid (www.westgrid.ca) and Compute Canada Calcul Canada (www.computecanada.ca).

\bibliographystyle{elsarticle-num-names}


\end{document}